\documentclass[12pt,leqno]{amsart}
\usepackage{amsmath,amssymb,array,latexsym}

\theoremstyle{plain}
\newtheorem{thm}{Theorem}[section]
\newtheorem{pro}[thm]{Proposition}
\newtheorem{lem}[thm]{Lemma}
\newtheorem{cor}[thm]{Corollary}
\theoremstyle{definition}
\newtheorem{defn}[thm]{Definition}
\newtheorem*{rem}{Remark}

\newcommand{\I}{\mathbb{I}}
\newcommand{\N}{\mathbb{N}}

\newcommand{\R}{\mathbb{R}}
\newcommand{\cA}{\mathcal{A}}
\newcommand{\cC}{\mathcal{C}}
\newcommand{\cB}{\mathcal{B}} 

\newcommand{\cF}{\mathcal{F}}
\newcommand{\cG}{\mathcal{G}}
\newcommand{\cH}{\mathcal{H}}
\newcommand{\cL}{\mathcal{L}}

\newcommand{\cS}{\mathcal{S}}
\newcommand{\cU}{\mathcal{U}}
\newcommand{\CB}{\text{CB}}
\newcommand{\CBS}{\text{CBS}}
\newcommand{\vp}{\varepsilon}
\newcommand{\Ep}{\text{Ep}}
\newcommand{\Id}{\text{Id}}
\newcommand{\In}{\text{In}}

\newcommand{\alphao}{{\alpha_0}}

\newcommand{\ub}{\overline{u}}
\newcommand{\vb}{\overline{v}}

\newcommand{\xb}{\overline{x}}
\newcommand{\yb}{\overline{y}}
\newcommand{\zb}{\overline{z}}

\newcommand{\Lim}{\text{Lim}}
\newcommand{\supp}{\text{supp}}
\newcommand{\sign}{\text{sign}}

\newcommand{\coo}{\text{c}_{00}}
\newcommand{\co}{\text{c}_{0}}

\newcommand{\infsupset}{\stackrel{\infty}{\supset}}
\newcommand{\infsubset}{\stackrel{\infty}{\subset}}

\title{How many operators do there exist on  a Banach space?}

\author{Th.~Schlumprecht}
\thanks{The research was supported by NSF}

\begin{document}

\begin{abstract}
 We present some results concerning the following  question:
Given an infinite dimensional Banach space $X$. Are there two
 normalized basic  sequences $(x_n)$ and $(y_n)$, with 
$\liminf_{n\to\infty} \|x_n-y_n\|>0,$ so
that the map  $x_n\mapsto y_n$
 extends to a linear bounded operator between the linear span
    $(x_n)$ and the linear span  $(y_n)$?
\end{abstract}

\maketitle
\markboth{TH. SCHLUMPRECHT}{OPERATORS ON BANACH SPACES}

\section{Introduction}\label{S:1}
Over the last years the following
 question caught the interest  of an increasing
number of researchers in the area of
 infinite dimensional Banach space theory and became
 for many of them one of the central problems.

\begin{enumerate}
\item[(Q1)]  Assume $X$ is an infinite dimensional
 Banach space. Does there exist a linear bounded operator $T$
 on $X$ which is not of the form $T=\lambda \Id +C$, where $\lambda$ is a scalar,
 $C$ is a compact operator on $X$, and $\Id$ denotes the identity on $X$?
\end{enumerate}

 Clearly, if $E\subset X$ is  a finite dimensional
 subspace and $S: E\to X$ is any linear operator, then $S$ can be extended 
 to a linear operator $T:X\to X$ of finite rank. Thus, there must be ``many'' finite  rank
operators on $X$, and thus  ``many'' elements in the closure
 of the finite rank operators, which is, assuming $X$ has the approximation property,
 the ideal of compact operators. Of course,
 the identity  must also be a linear bounded operator. Therefore question
 (Q1) asks wheather or not
   there are always some other linear bounded operators defined on
  a Banach space,
 beyond the class of operators which must exist by elementary reasons.

 Note that a counterexample to (Q1) would be the first example of a
 Banach space $X$ for which the invariant subspace problem
has a positive answer: Does every operator $T$ on $X$ have a non trivial
 subspace $Y$ ($Y\not=\{0\}$  and $Y\not=X$) so that $T (Y)\subset Y$?
 Indeed, due to a result of V.~I.~Lomonosov \cite{Lo}, any operator
on a Banach space which commutes with  a compact operator $C\not= 0$ must
have a none trivial invariant  subspace.
 A counterexample to (Q1) might also be a prime candidate for a Banach space
 on which all Lipschitz function admit a point of differentiability \cite{Li}.

Recall (\cite{GM1} and \cite{GM2}) that an infinite dimensional  Banach space $X$
is called a {\em hereditarily indecomposable space} (HI) if
 every operator $T:Y\to X$, with $Y$ being a subspace of $X$ is of the form
$\In_{(Y,X)}+S$, where $\lambda$ is a scalar, $\In_{(Y,X)}$ denotes the inclusion map from $Y$ into $X$ and
$S$ is a {\em strictly singular operator}, i.e. an operator which on 
no infinite dimensional subspace   is an isomorphism into $X$.

The first known example of an HI space, we denote it by $GM$,  was constructed by T.~W.~Gowers and
 B.~Maurey  in  1993 \cite{GM1}. Since then HI spaces with additional properties
were constructed (e.g. \cite{AD}, \cite{F}, and  \cite{ OS}).

 Although a counterexample to (Q1) does not need to be an HI
space, these spaces are nevertheless natural candidates for
 counterexamples. In \cite{AS} it was shown that the space constructed
in \cite{GM1} admits an operator (defined on all of $GM$) which
is strictly singular but not compact. Already earlier it was shown in \cite{Go1} that
 such an operator can be constructed on a subspace of $GM$.

 Given the seemingly easy task to write down on a concrete space a nontrivial
operator (i.e. not of the form $\lambda+compact$) the proof in \cite{AS} is
quite involved. It would be interesting, but probably technically even harder,
 to establish the existence of non trivial operators defined on  other known  HI spaces.
 Since it seems to be hard to answer (Q1) already for specific spaces, the tools to
  give an answer (at least a positive one) in the general case   are probably
  not developed yet.  It also seems that (Q1) is the type of question which will not
 be answered directly, but by solving first some other, more structure theoretical questions.

Following more along the line of the  concept of HI spaces we  turn therefore to the
 following ``easier question'' (for a positive answer) and ask whether or not it is
always possible to define a non trivial operator on a subspace.

\begin{enumerate}
\item[(Q2)] Assume $X$ is an infinite dimensional Banach space. Is there  a closed subspace $Y$
 of $X$  and an operator $T:Y\to X$, so that $T$ is not of the form
 $T=\lambda \In_{(Y,X)}+ C$, where $\lambda$ is a scalar and  $C:Y\to X$ is a compact operator?
\end{enumerate}

It is easy to see, that (Q2) can be equivalently reformulated as follows

\begin{enumerate}
\item[(Q3)] Assume that $X$ is an infinite dimensional Banach space. Does there
exist two normalized basic sequences $(x_n)$ and $(y_n)$ so that $(x_n)$ dominates $(y_n)$
 and so that $\inf_{n\in\N}\|x_n-y_n\|>0$?
\end{enumerate}

We want to go one step further and formulate a more structure theoretical approach to our problem.
We call a  
Banach space $X$ with a basis $(e_i)$ {\em a space of Class 1} if:
\begin{enumerate}
\item[(C1)] Every block basis of $(e_i)$ has a subsequence
 which is equivalent to some subsequence of $(e_i)$.
\end{enumerate}

\begin{rem}
Note that the spaces $\ell_p$, $1\le p<\infty$, and $c_0$ are clearly
 (C1). Actually in  $\ell_p$, $1\le p<\infty$, and $c_0$ all blockbasis
 are equivalent to each other, a property which characterizes the unit bases
 of  $\ell_p$, $1\le p<\infty$, and $c_0$ \cite{Z}.
 Moreover, in \cite{CJT} it was shown that Tsirelson's space $T$
 (as described in \cite{FJ}) as well as its dual $T^*$
 (the original {\em Tsirelson space}  defined in \cite{T}) are spaces of Class 1. 
 Recently  the result of \cite{CJT} was generalized to all
{\em finitely mixed Tsirelson spaces} in \cite{LT}. 
  The reader unfamiliar with the usual notations
  and concepts  (like  mixed Tsirelson spaces) is referred   to the last paragraphs of this  section,
where all the notions  of this paper are introduced.
 
Also note that every space of Class 1 must contain 
 an unconditional basic sequence, an observation which follows immediately
 from the result in \cite{Go1} which says that every infinite dimensional
 Banach space must either contain an infinite dimensional subspace with
unconditional basic sequence or it must contain an infinite dimensional
 subspace which is HI. Of course, a space of Class 1 cannot contain
an infinite dimensional  subspace which is HI.
 
We do not know of any elementary proof, i.e. a proof which does not use the above cited result
 of \cite{Go1}, of the fact that a space of Class 1 must contain  an unconditional basic sequence.
\end{rem}

It  seems that  until the early nineties all known Banach spaces
 had subspaces which were  of Class 1.
 Then, in  1991, the author \cite{Sch} of this work constructed a space, nowadays denoted by  $S$,
 which fails (C1) in an extreme way: 
In all infinite dimensional subspaces, spanned by a block, one is not only able to find
 two normalized blocks  $(x_n)$ and $(y_n)$ 
 having no subsequences which are equivalent to each other.
 But one can even find in each  subspace
 two normalized blocks  $(x_n)$ and $(y_n)$, for which the map $x_n\mapsto y_n$ extends
to a bounded linear  and strictly singular operator between the span of $(x_n)$ and  the span of $(y_n)$.

 Let us therefore define the following second  class of Banach spaces.
 Let us call a Banach space $X$ with a basis $(e_i)$  a space of Class 2 if:

\begin{enumerate} 
\item[(C2)] Each block basis $(z_n)$ has two further block bases $(x_n)$
 and $(y_n)$  so that the map $x_n\mapsto y_n$   extends to a bounded linear and strictly
 singular operator  between the span of   $(x_n)$ and the span of $(y_n)$.
\end{enumerate}

The main purpose of this paper is to establish  criteria  for a Banach space to be
 of Class 2 and thereby  address the following  problem.

\begin{enumerate}
\item[(Q4)] Does every infinite dimensional Banach space contain a Banach space which is either
 of Class 1 or of Class 2?
\end{enumerate}

Question (Q4) might seem at first sight   somewhat daring, let us therefore motivate it.  
If one believes (Q3) to be true, one could argue that experience
 shows us  that positive results in Banach space theory are often
derived from dichotomy principles. 
 Therefore (C2) could be    a   candidate  for the second alternative in a dichotomy
 result in which (C1) is the first alternative.

Secondly, if one believes Question (Q3) to have a counterexample, one might 
 try to look first for a  counterexample of (Q4).   Question (Q4) could have, contrary to
 Question (Q3),  a counterexample with an unconditional basis, which  probably
 will be easier to define.  Then, starting with a counterexample to (Q4), one
 might ask for a modification of it to obtain a counterexample to question (Q3).
Similarly, the key argument toward defining the first known HI space,
 was the  definition of a  space of Class 2.

Finally note that a negative answer of (Q3) is equivalent to the statement
 that the following class of Banach spaces is not empty.
A Banach space $X$ with a basis $(e_i)$ is said to be {\em a space of Class 3} if: 
\begin{enumerate}
\item[(C3)]  For any two block bases $(x_n)$ and $(z_n)$ of $(e_i)$, with 
   $\inf_{n\in\N}\|x_n-z_n\|>0$, neither the map $x_n\mapsto z_n$ nor
 the map $z_n\mapsto x_n$ extends to a linear bounded operator between 
  the spans of  $(x_n)$ and $(z_n)$.
\end{enumerate}

In this paper we are interested in formulating criteria which imply that
 a given Banach space $X$ is of Class 2. Therefore we want to find sufficient conditions
 for $X$ to contain two seminormalized basic sequences $(x_n)$ and $(y_n)$ for which
 the map $x_n\mapsto y_n$ extends to a bounded linear and strictly singular
 operator.

In  Section \ref{S:2} we discuss the following problem: Given
two normalized basic sequences $(x_n)$ and $(y_n)$, which conditions should we impose
 on the spreading models of $(x_n)$ and $(y_n)$, in order to insure
the existence of subsequences $(\tilde x_n)$ and $(\tilde y_n)$ of
$(x_n)$ and $(y_n)$ respectively, so that $\tilde x_n\mapsto\tilde y_n$
 extends to a strictly singular and bounded linear operator?
 The properties of spreading models can be  quite different form the properties
 of the underlying  generating sequences. For example, the spreading model of the basis of
 Schreier's space (using the notations of Definition \ref{D:1.4} this is the space
 $S(\cS_1)$, with $\cS_1=\{E\in[\N]^{<\infty}: \#E\le \min E\}$)
is isometrically equivalent to the $\ell_1$-unit vector basis. On the other
hand Schreier's space is hereditarily $\co$.
Therefore we expect to need to impose rather strong conditions on the spreading
 models of $(x_n)$ and $(y_n)$  in order to conclude that for some
 subsequences  $(m_k),\quad(n_k)$ and $\N$ it follows that
 $x_{m_k}\mapsto y_{n_k}$ extends to bounded linear  and strictly singular map
$[x_{m_k}:k\in\N]\to[y_{n_k}:k\in\N]$.

The main result of  Section \ref{S:2} is the following answer to our question.

\begin{thm}\label{T:2.0}
Let$(x_n)$ and $(y_n)$ be two normalized weakly null sequences
having spreading models $(e_n)$ and $(f_n)$ respectively. 

Assume  that $(e_n)$ is not equivalent to the
 $c_0$-unit vector basis and that the following condition holds.
\begin{align} &\text{There is a sequence  $(\delta_n)$ of positive numbers decreasing
 to zero so that}\qquad\\
&\qquad\qquad\qquad\Big\|\sum_{i=1}^\infty a_if_i\Big\|\le
\max_{\substack{n\in\N \\  i_1<i_2<\ldots i_n}}\delta_n\Big\|\sum_{j=1}^n a_{i_j}e_i\Big\|.\notag
\end{align}
Then there are subsequences $(\tilde x_n)$ and $(\tilde y_n)$
of $(x_n)$ and $(y_n)$ respectively, so that   $\tilde x_n\mapsto\tilde y_n$ extends
 to a strictly singular and bounded linear operator.
\end{thm}

 In order to formulate our second result, we need to recall
 the Cantor Bendixson index of subsets of $[\N]^{<\infty}$, the set of 
 finite subsets of $\N$.

 A set $\cA\subset[\N]^{<\infty}$ is called
{\em hereditary} if for all $A\in\cA$ and all $B\subset A$ it 
follows that $B\in \cA$. We call $\cA$ {\em spreading} if
   $A=\{a_1,a_2,\ldots a_k\}\in \cA$, with $a_1<a_2<\ldots <a_k$, and
if $b_1<b_2<\ldots b_k$ are in $\N$, so that $b_i\ge a_i$
for $i=1,\ldots k$, then $\{b_1,\ldots b_k\}\in\cA$. We
 always consider  the topology
 of pointwise convergence on $[\N]^{<\infty}$(identifying $A\subset\N$ with $\chi_A$,
 the characteristic function of $A$).

For an hereditary and compact $\cA\subset[\N]^{<\infty}$ we define
\begin{align}\label{E:1.1}
\cA^{(1)}&=\{ A\in\cA| A\text{ is accumulation point of }\cA\}\\
&=\{ A\in\cA|\exists N\subset\N, N \text{  infinite }\forall n\in N\quad A\cup\{n\}\in\cA\}.\notag
\end{align}
Since $\cA$ is compact in $[\N]^{<\infty}$ each $A\in\cA$ is subset of a maximal element of $\cA$
and therefore we conclude that
\begin{equation}\label{E:1.2}
\cA^{(1)}\subsetneq \cA,\text{ if }\cA\not=\emptyset.
\end{equation}
Secondly it follows that $\cA^{(1)}$ is also compact and hereditary. We can therefore
 define $ \cA^{(\alpha)}$ for each $\alpha<\omega_1$ by transfinite induction.
 We let $ \cA^{(0)}=\cA$ and assuming we defined $\cA^{(\beta)}$ for all $\beta<\alpha$ we
 put
 \begin{equation}\label{E:1.3}
  \cA^{(\alpha)}=\bigl(\cA^{(\beta)}\bigr)^{(1)}\text{ if }\alpha=\beta+1,\text{ and }
   \cA^{(\alpha)}=\bigcap_{\beta<\alpha} \cA^{(\beta)}\text{ if }\alpha=\sup_{\beta<\alpha} \beta.
\end{equation}
Since $[\N]^{<\infty}$ is countable it follows from (\ref{E:1.2}) that there is an $\alpha<\omega_1$ for
 which $\cA^{(\alpha)}$ is empty and we define the Cantor Bendixson index of $\cA$ by
 \begin{equation}\label{E:1.5}
\CB(\cA)=\min\{\alpha<\omega_1: \cA^{(\alpha)}=\emptyset\}.
\end{equation}
We note that $\CB(\cA)$ is a successor ordinal. Indeed, assume that $\alpha<\omega_1$ is a limit
 ordinal and that $\cA^{(\beta)}\not=\emptyset$ for all $\beta<\alpha$. Then it follows from
the fact that  $\cA^{(\beta)}$ is hereditary that $\emptyset\in
\cA^{(\beta)}$ for all $\beta<\alpha$, and, thus, by (\ref{E:1.3}) that
 $\emptyset\in\cA^{(\alpha)}$ which implies that $\CB(\cA)>\alpha$.

The {\em strong Cantor Bendixson index}  of $\cA$ is defined as follows (cf. \cite{AMT}):
 \begin{equation}\label{E:1.6}
\CBS(\cA)=\sup_{N\subset \N, \#N=\infty}\inf_{M\subset N, \#M=\infty} \CB\bigl(\cA\cap[M]^{<\infty}\bigr)
\end{equation}

\begin{rem} Note that $\CBS(\cA)$ could be a limit ordinal. Indeed, let $(N_n)_{n\in\N}$ be a sequence of
infinite pairwise disjoint subsets of $\N$ whose union is $\N$.
 Take $\cA=\bigcup_{n\in\N} \cA_n$, with $\cA:=\{A\subset N_n\text{ finite }: \#A=n\}$.
Then it follows that $\CBS(\cA)=\omega$.

 On the other hand it is clear that for a successor ordinal $\gamma$  and
 an hereditary $\cA\subset[\N]^{<\infty}$ we have that
 \begin{equation*}
\CBS(\cA)\ge \gamma\iff \exists N\subset\N,\#N=\infty\forall M\subset N,\#M=\infty\quad
 \CB(\cA\cap[M]^{<\infty})\ge \gamma.
\end{equation*}

\end{rem}

Using the Cantor Bendixson index  we can characterize spaces  which have
subspaces of Class 1 and quantify the property of not having a subspace of Class 1. More
 generally, we will use the Cantor Bendixson index to measure ``how far away''
 two basic sequences are to each other, up to passing to subsequences''.
\begin{defn}\label{D:1.1a}
Assume that $\xb=(x_n)$ and $\zb=(y_n)$ are two seminormalized basic sequences.
 For $c\ge 1$ we define
\begin{equation}\label{1.1a.1}
\cC(\xb,\zb,c)=\{A\in[\N]^{<\infty}: \exists B\in [\N]^{<\infty},\#B=\#A, \quad(x_n)_{n\in A}\sim_c (z_n)_{n\in B}\},
\end{equation}
by ``$(x_n)_{n\in A}\sim_c (z_n)_{n\in B}$'' we mean that $(x_n)_{n\in A}$ and $(z_n)_{n\in B}$
 are $c$-equivalent (see end of section), where on $A$ and $B$ we consider the order given by $\N$.
For $\gamma<\omega_1$ let (put $\inf\emptyset=\infty$)
\begin{align}\label{1.1a.2}
c(\xb,\zb,\gamma)&=\inf\bigl\{ c\ge 1: \CBS\bigl(\cC(\xb,\zb,c)\bigr)\ge\gamma+1\bigl\}
   \text{ (with $\inf\emptyset=\infty$)}\\
                 &= \inf\bigl\{ c\ge 1:\exists N\infsubset \N\,\forall M\infsubset N\quad
                    \CB\bigl(\cC(\xb,\zb,c)\cap[M]^{<\infty}\bigr)\ge\gamma+1\bigl\}
                          \notag\\
\label{1.1a.3}
\gamma_0(\xb,\zb)&=\sup\{ \gamma<\omega_1:c(\xb,\zb,\gamma)<\infty\}.
\end{align}
In the  case that $\zb$ is the unit basis of $\ell_1$ we define
 $\cB(\xb,c)=\cC(\xb,\zb,c)$, $b(\xb,\gamma)=c(\xb,\zb,\gamma)$, and
 $\beta_0(\xb)=\gamma_0(\xb,\zb)$ and for a Banach space $X$ and $\beta<\omega_1$ we put
\begin{align}\label{1.1a.4}
 B(\beta)=B(\beta, X)&=\sup\{b(\zb,\beta) : \zb=(z_n)\subset B_X \quad \text{semi normalized, weakly null}\}\\
   &=\sup\Biggl\{ b>0:
 \begin{matrix} & \exists (z_n)\subset B_X \quad \text{semi normalized and weakly null}\\
                & \CBS(\cB(\zb,b))\}\ge \beta+1
\end{matrix}\Biggr\}.\notag
\end{align}
\begin{align}\label{1.1a.5}
\beta_0(X)&=\sup\{ \beta\!<\omega_1 : B(\beta)\!>\!0\}\\
&=\sup\{ \beta_0(\zb):\zb=(z_n)\subset B_X \quad \text{semi normalized and weakly null}\}\notag
\end{align}
It is  clear that $\beta_0(X)$ must be a limit ordinal.
\end{defn}

\begin{pro}\label{P:1.1}(see proof after Corollary \ref{C:5.6}  in Section \ref{S:4}).
Assume $X$ is a Banach space with  a basis $(x_n)$.
If $(z_n)$ is another basic sequence for which
 $\gamma_0((z_n),(x_n))=\omega_1$ then
 a subsequence of $(z_n)$ is
 ismorphically equivalent to a subsequence of $(x_n)$ or
 to the spreading model of a subsequence of $(x_n)$.

 Moreover, $X$ is of Class 1 if and only if
 for all blockbases $(z_n)$ of $(x_n)$ it follows that
$\gamma_0((z_n),(x_n))=\gamma_0((x_n),(z_n))=\omega_1$.
\end{pro}
In Section  \ref{S:6} we will prove the following result,
  implying a condition for a space to be in Class 2.

\begin{thm}\label{T:1.2}
Let $X$ be a Banach space with a basis $(e_i)$ not containing $c_0$. Assume that there is an
 ordinal $\gamma\in[0,\beta_0(X)]$  so that:
\begin{enumerate}
\item[a)] There is a semi normalized block basis $\yb=(y_n)\subset B_X$ with
 $\beta_0(\yb)\le\gamma$.
\item[b)] $\inf_{\beta<\gamma} B(\beta,X)>0$.
\end{enumerate}
Then there is a seminormalized block basis $(x_n)$ in $X$
 and a subsequence $(\tilde y_n)$ of $(y_n)$ so that
the map $x_n\mapsto \tilde y_n$ extends to a strictly singular operator
 $T: [x_n:n\in\N]\to [\tilde y_n:n\in\N]$.
\end{thm}

\begin{cor}\label{C:1.2a}
Assume that $Z$ is a reflexive Banach space with a basis and each block subspace
 of $Z$ satisfies the condition of Theorem \ref{T:1.2}.

Then $Z$ is of Class 2.
\end{cor}

\begin{rem}  Note that in the statement of Theorem \ref{T:1.2}
either $\gamma=\beta_0(X)$ in which case it is clear that (a) is satisfied for
all block bases. Or $\gamma<\beta_0(X)$ in which case it is clear that
 (b) is satisfied.

Also note that the assumption in Theorem \ref{T:1.2} says the following:

On the one hand there is a block basis $\yb$ which is ``far away from the $\ell_1$-
unit vector basis, in the sense that given any $\vp>0$ there is an $\alpha<\gamma$ so that
  $b(\yb,\alpha)<\vp$, i.e. on sets of strong Cantor Bendixson index $\alpha$ 
 the best equivalence constant between the $\ell_1$ basis and $\yb$ is at most $\vp$.

On the other hand if we let $c=\inf_{\beta<\gamma} B(\beta,X)$, we can choose
 for any $\alpha<\gamma$ and any $\vp>0$
a block sequence $\zb$  for which    $b(\zb,\alpha)>c-\vp$.

In other words the  assumption in Theorem \ref{T:1.2} states that we require 
 the existence of blockbases which have ``different $\ell_1$-behavior''.

But if one wants to attack Problem (Q4) one could start by assuming
that our given  Banach space $X$ has no  infinite dimensional subspace of
Class 1 and therefore by Proposition \ref{P:1.1}
 every block basis $(y_n)$ must admit a further block basis
 $(z_n)$ so that 
$\gamma_0((y_n),(z_n))<\omega_1$ 
or $\gamma_0((z_n),(y_n))<\omega_1$. So, instead of going the quite ``uneconomical route''
 (as in Theorem  \ref{T:1.2}) and  comparing two blockbases to the $\ell_1$-unit vector basis
 one should compare them to each other directly.

Here exactly lies one major hurdle which one has to overcome before
 reaching further results. The proof of Theorem \ref{T:2.0} as well as the proof of
 Theorem \ref{T:1.2} make use of certain ``partial unconditionality results''
 namely a result by E.~Odell \cite{O2} (see Theorem \ref{T:2.2a})
 and a result by S.~A.~Argyros, S.~Mercourakis and A.~Tsarpalias
 \cite{AMT} (see Lemma \ref{L:6.1} in Section \ref{S:6}). Both results
 give partial answers to the question, under which conditions
 a family  of finite rank projections on certain subsets of a basis
 are uniformly bounded.
 
In order to make advances on (Q4) we would need to address the
following question which asks for generalizations of the results
 in \cite{O2} and \cite{AMT}.

\begin{enumerate}
\item[(Q5)] Assume that $X$ is a Banach space with a basis $(e_n)$, $\gamma\in[0,\omega_1)$,
   and  assume that for  all
 normalized block bases $(x_n)$ of $(e_n)$ it follows that
 $c((x_n),(e_n) ,\gamma)<\infty$ 
 $c((e_n) ,(x_n),\gamma)<\infty$.
  
 Secondly, consider for a blockbasis $(x_n)$ and an $r> 0$ the set
 
$$\cU((x_n),r)=\big\{ A\in[\N]^{<\infty}: \|P_{[x_n:n\in A]}\|\le r\big\},$$
where $P_{[x_n:n\in A]}:[x_n:n\in \N]\to [x_n:n\in A]$ is the defined to be
the usual projection onto $[x_n:n\in A]$.

Now, does it follow that any block basis $(x_n)$ has a subsequence 
 $(z_n)$ for which $\CBS(\cU((z_n),r))\ge \gamma$ for some $r>0$?
 
\end{enumerate}

The result cited from \cite{O2}  gives a positive answer 
 if $\gamma=\omega$ and the result cited from \cite{AMT} 
  gives a positive answer  if we replaced in the assumption 
 $c((e_n) ,(x_n),\gamma)$  by $b((x_n),\gamma)$.
                                                                   
\end{rem}

 To be able to prove Theorem \ref{T:1.2} we will introduce in Section \ref{S:3} a  family
 $(\cF_\alpha)_{\alpha<\omega_1}$ of subsets of $[\N]{<\infty}$,
having the property that for each $\alpha<\omega_1$ the strong
Cantor Bendixson index is $\alpha+1$ (see Corollary\ref{C:5.6}).
 The family $(\cF_\alpha)$ is defined in a similar fashion as
the {\em Schreier families} $(\cS_\alpha)_{\alpha<\omega_1}$
 \cite{AA}, with one crucial difference:
 one obtains $\cF_{\alpha+1}$ from  $\cF_{\alpha}$
 by  adding to each element of $A\in\cF_\alpha$ at most one
new element, therefore the step from $\cF_\alpha$
 to  $\cF_{\alpha+1}$ is much smaller than
 the step from  $\cS_\alpha$ to  $\cS_{\alpha+1}$. 

 Feeling that this family
 $(\cF_\alpha)$  could be an important tool to analyse the
combinatorics of blockbases we extensively discuss the properties
 of this family in the Sections \ref{S:3}, \ref{S:4} and \ref{S:5}, in particular we prove 
 a result which, roughly speaking, says that
 any hereditary  set $\cA\subset[\N]^{<\infty}$ restricted to some appropriate 
$M\infsubset\N$ is  equal to ``a
 version  of $\cF_\alpha$'' (see Theorem \ref{T:4.1}).
Based on this result we will be able to prove Theorem \ref{T:1.2} in Section \ref{S:6}.

We  will need the following notations and conventions.

For simplicity all our  Banach spaces are considered to
 be over the field $\R$. If $X$ is a Banach space $B_X$ and $S_X$ denotes
the unit ball and the unit sphere respectively.

$\coo$ denotes the  vectorspace of sequences in $\R$ which eventually vanish and
 $(e_i)$ denotes the usual vector basis of $\coo$ and when we consider a Banach space $X$
with a normalized basis, we think of $X$ being the completion of some norm on $\coo$
 with $(e_i)$ being that basis.
 For  $x=\sum x_i e_i\in \coo$
 we let $\supp(x)=\{i\in\N: x_i\not=0\}$ be the {\em support of $x$} and if $E\subset \N$  we let
$E(x)=\sum_{i\in E} x_i e_i$.

We say that a basic sequence $(x_n)$ dominates another basic sequence 
$(y_n)$ if the map $x_n\mapsto y_n$ extends to a  linear bounded map between
the span of $(x_n)$ and the span of $(y_n)$, i.e. if the there is a $c>0$ so that for all
 $(a_i)\in\coo$
 $\|\sum_{i=1}^n a_i y_i\|\le c\|\sum_{i=1}^n a_i x_i\|$.

We say that  $(x_n)$ and  $(y_n)$ are $c$-equivalent, $c\ge 1$ and we write $(x_n)\sim_c(y_n)$
 if for any $(a_i)\in\coo$ it follows that
 $
  \frac1c \|\sum_{i=1}^n a_i x_i\|\le  \|\sum_{i=1}^n a_i y_i\|\le c\|\sum_{i=1}^n a_i x_i\|$.

The closed linear span of a subsequence $(x_n)_{n\in\N}$ of a Banach space is denoted by
 $[x_n\!:\!n\!\in\!\N]$.

 Let $E$ be a Banach space with a 1-spreading basis $(e_i)$, i.e. 
 \begin{equation*}
\Bigl\|\sum_{i=1}^m a_i e_i\Bigr\|=\Bigl\|\sum_{i=1}^m a_i e_{n_i}\Bigr\|,
\text{ whenever $m\in\N$, $(a_i)_{i=1}^m\subset\R$ and $n_1<n_2<\ldots n_m$ are in $\N$},
\end{equation*}
and let $(x_n)$ be a seminormalized basic sequence. We say that  $(e_i)$ is the spreading model of 
  $(x_n)$ if for any $k\in\N$ and any $(a_i)_{i=1}^k\subset \R$ it follows that
 \begin{equation*}
\lim_{n_1\to\infty}\lim_{n_2\to\infty}\ldots\lim_{n_k\to\infty}
\Bigl\| \sum_{i=1}^k a_i x_{n_i}\Bigr\|=\Bigl\|\sum_{i=1}^k a_ie_i\Bigr\|.
\end{equation*}
Recall \cite{BL} that any seminormalized basic sequence has a subsequence with
 spreading model and that, using our notations in Definition \ref{D:1.1a},
 a seminormalized sequence $(x_n)$ has  a subsequence whose spreading model
 is isomorphic to $(e_i)$ if and only if  $c((x_n),(e_i),\omega)<\infty)$.

For a set $M$ the set of all finite subsets of $M$ is denoted by $[M]^{<\infty}$
 and the set of all infinite subsets is denoted by $[M]^{\infty}$.
 We write $M\infsubset N$ if M is an infinite subset of $N$.

\begin{defn}\label{D:1.4} (Mixed Schreier and Tsirelson spaces)

For an  hereditary, spreading and
 compact   subset $\cF$ of $[\N]^{<\infty}$
 containing all singletons of  $[\N]^{<\infty}$
 we define {\em the $\cF$-Schreier space} $S(\cF)$
    to be the completion of $\coo$ under the norm defined by

\begin{equation}
  \|x\|=\sup_{E\in\cF} \sum_{i\in E} |x_i|, \text{ whenever }x=(x_i)\in\coo.
\tag{$S(\cF)$} 
\end{equation}
 Let $(\cF_n)_{n=1}^\infty $ be an increasing sequence of  hereditary, spreading and
 compact   subsets of $[\N]^{<\infty}$ and  let $(\Theta_n)\subset(0,1]$
be a non increasing sequence.
The {\em mixed Schreier space } $S((\Theta_n),(\cF_n))$  is the completion of $\coo$ 
under the norm defined by
\begin{equation}
\|x\|=\sup_{n\in\N}\sup_{E\in\cF_n}\Theta_n\sum_{i\in E} |x_i|, \text{ whenever }x=(x_i)\in\coo.
\tag{$S((\Theta_n),(\cF_n))$}
\end{equation}

If in addition $(\ell_n)_{n\in\N}$ is a sequence in $\N$  increasing to  $\infty$,
 we define  the {\em  mixed Schreier space with the additional admissibility condition
 given by  $(\ell_n)_{n\in\N}$} to be the space
  $S((\Theta_n),(\cF_n),(\ell_n))$  with the following norm

\begin{equation}
\|x\|=\sup_{n\in\N}\sup_{\ell_n\le\min E, E\in\cF_n}\Theta_n\sum_{i\in E} |x_i|,
\quad \text{ whenever }x=(x_i)\in\coo.\tag{$S((\Theta_n),(\cF_n),(\ell_n))$}
\end{equation}

If $\cF$ of $[\N]^{<\infty}$ is hereditary, spreading and compact, and if
 $E_1<E_2<\ldots E_\ell$ are in $\cF$, we say that
 $(E_i)_{i=1}^\ell$ is $\cF$-{\em admissible} if
 $\{\min E_i: i=1,2\ldots \ell\}\in\cF$.
The {\em mixed Tsirelson space } $T((\Theta_n),(\cF_n))$  is the
completion of $\coo$
under the norm  which is implicitly defined by
\begin{equation}
\|x\|=\max\left\{\sup_{n\in\N}|x_n|,\sup_{n\in\N}\sup_{(E_i)_{i=1}^n\text{ $\cF_n$-adm.}} 
\Theta_n\sum_{i=1}^n \|E_i(x)\|\right\} \text{ for }x\!=\!(x_i)\!\in\!\coo
\tag{$T((\Theta_n),(\cF_n))$}
\end{equation}
A {\em finitely mixed Tsirelson space } is a space  $T((\Theta_n),(\cF_n))$
 where $(\Theta_n)\subset(0,1]$ and $(\cF_n)$ are finite sequences.
\end{defn}

\section{Conditions on the spreading models implying the existence
 of nontrivial operators}\label{S:2}

\begin{defn}\label{D:2.1} For two normalized basic sequences $(x_n)$ and $(y_n)$
 and an $\vp>0$ we define
\begin{equation*}
\Delta_{((x_n),(y_n))}(\vp)=\sup
\left\{\Big\|\sum a_i y_i\Big\|:
(a_i)\in\coo,\, \max_{i\in\N}|a_i|\le\vp,\text{ and }  \Big\|\sum a_i x_i\Big\|\le 1 \right\}.
\end{equation*}
We say that $(x_n)$ {\em strongly dominates } $(y_n)$ if
\begin{align}\label{E:2.1.1}
  &\lim_{n\to\infty}\inf_{A\subset \N,\#A=n}\|\sum_{i\in A} x_i\|=\infty\text{ and }\\
 \label{E:2.1.2} & \lim_{\vp\searrow 0}\Delta_{((x_n),(y_n))}(\vp) = 0.
\end{align}
\end{defn}

\begin{rem} Assume $(x_n)$ and $(y_n)$ are two normalized basic sequences.
\begin{enumerate}
\item[a)]  If $\Delta_{((x_n),(y_n))}(\vp)<\infty$ for some $\vp>0$
 (and thus for all $\vp>0$) then $(x_n)$ dominates $(y_n)$,
i.e. the map $x_n\mapsto y_n$ extends to a bounded linear operator.
\item[b)] Note that in the case that $(x_n)$ is subsymmetric condition (\ref{E:2.1.1})
 means that $(x_n)$ is not equivalent to the $\co$ unit vector basis.
 \item[c)] No  normalized basic sequence strongly dominates the unit vector basis in $\ell_1$.
\end{enumerate}
\end{rem}

 We will need the following result
 by E.~Odell on {\em Schreier unconditionality}.
\begin{thm}\label{T:2.2a}\cite{O2} Let $(x_n)$ be a normalized weakly null sequence in a Banach space and $\eta>0$.
 Then $(x_n)$ contains a subsequence $(\tilde x_n)$ so that

 \begin{equation}\label{E:2.2a.1}
  \forall (\alpha_i)\in\coo \forall F\subset \N, \min F\ge \#F:\quad
  \left\|\sum_{i\in F} \alpha_i \tilde x_i\right\|\le
 (2+\eta)\left\|\sum_{i=1}^\infty \alpha_i \tilde x_i\right\|.\notag
\end{equation}
If $(x_n)$ is a bimonotone sequence above factor $(2+\eta)$ could be replaced
by $(1+\eta)$.
\end{thm}

\begin{rem} Assume that $(x_n)$ and $(y_n)$ are normalized basic sequence,
 that $(x_n)$ strongly dominates $(y_n)$ and that
 either $(x_n)$ is weakly null  or 1-subsymmetric. 

By passing to subsequences
 (of $(x_n)$ and $(y_n)$ simultaneously)
we can assume that
 $(x_n)$ has a spreading model $(e_i)$ and that
the conclusion of Theorem \ref{T:2.2a} holds, i.e. that
\begin{equation}\label{E:2.2a.2}
\|\sum_{i\in F} \alpha_i e_i\|\sim_2\|\sum_{i\in F} \alpha_i x_i\|
\le 3\|\sum_{i=1}^\infty \alpha_i x_i\|
\end{equation}
for all $(\alpha_i)\in\coo$ and  finite $F\subset\N$, with $\min F\ge \#F$.

Since $(x_n)$  satisfies (\ref{E:2.1.1}) $(e_n)$ cannot be equivalent to the
 $c_0$-unit basis and we deduce therefore that
\begin{equation}\label{E:2.2a.3}
c_n=\|\sum_{i=1}^n e_i\|\nearrow\infty.
\end{equation}
We define now the following basic sequence $(\tilde y)$ which dominates $(y_n)$
 by
\begin{equation}\label{E:2.2a.4}
 \|\sum_{i=1}^\infty \alpha_i \tilde y_i\|=
 \|\sum_{i=1}^\infty \alpha_i y_i\|+
\max_n\frac1{\sqrt{c_n}}\max_{k_1<\ldots k_n}\|\sum_{i=1}^n\alpha_{k_i} e_i\|
\text{ whenever $(\alpha_i)\in\coo$}.
\end{equation}
\end{rem}
\begin{pro}\label{P:2.2b}
Using the assumptions and notations as in above
Remark  it follows that
\begin{enumerate}
\item[a)] $(x_n)$ strongly dominates $(\tilde y_n)$.
\item[b)] For $\eta>0$ it follows that
 $$\ell_\eta=\bigl\{\#\{i\in\N:|\alpha_i|\ge \eta\}:
  (\alpha_i)\in\coo,\|\sum\alpha_i\tilde y_i\|\le 1\bigr\}<\infty.$$
\end{enumerate}
\end{pro}
\begin{proof}
In order to show (a) let $\eta>0$. Since $(x_n)$ strongly dominates $(y_n)$
 we can choose an $\vp'>0$ so that $\Delta_{((x_n),(y_n))}(\vp')<\eta/2$. 
 Secondly we choose
 $n_0\in\N$ so that $1/\sqrt{c_{n_0}}<\eta/12$ and we let
$\vp=\min(\vp',\eta/2n_0)$.

If $(\alpha_i)\in\coo$ with $\|\sum\alpha_i x_i\|\le1$ and $\max|\alpha_i|\le \vp$
we deduce for $n\in\N$, with $n\le n_0$, that
 $$\frac1{\sqrt c_n}\max_{n\le k_1<\ldots k_n}
 \|\sum_{i=1}^n \alpha_i e_i\|
  \le {n\vp}\le\eta/2.$$
If  $n\in\N$, with $n\ge n_0$, we  deduce from (\ref{E:2.2a.2}) that
$$\frac1{\sqrt{c_n}}\max_{n\le k_1<\ldots k_n}\left\|\sum_{i=1}^n \alpha_{k_i} e_i\right\|\le
  \frac2{\sqrt{c_n}}\max_{n\le k_1<\ldots k_n}\left\|\sum_{i=1}^n \alpha_{k_i} x_i\right\|\le
  \frac6{\sqrt{c_n}}\left\|\sum_{i=1}^\infty \alpha_i x_i\right\|\le\eta/2,$$
which together with the choice of $\vp'$ implies
 that 
 $\|\sum\alpha_i \tilde y_i\|\le eta$ and finishes the proof of
the claim in (a).

To show (b), let $(\alpha_i)\in\coo$ so that
 $\|\sum \alpha_i\tilde y_i\|\le 1$
 and let $\eta>0$. We put $I=\{i\in\N:\|\alpha_i\|>\eta\}$
  and let $I'$ the set of the $\lceil\#I/2\rceil$ largest elements
 of $I$. Note that $\#I'\le \min I$ which implies by our
definition of  $(\tilde y_n)$
 and the fact that $(e_i)$ suppression 1-unconditional (thus 
 2-unconditional)
 that
\begin{align*}
1\ge\left\|\sum \alpha_i\tilde y_i\right\|\ge 
 \frac1{\sqrt{c_{\#I'}}}\left\|\sum_{i\in I'} \alpha_i e_i\right\|
\ge  \frac1{2\sqrt{c_{\#I'}}}\left\|\sum_{i\in I'} \eta e_i\right\|\ge 
 \eta\sqrt{c_{\#I'}}/2,
\end{align*}
which implies that $\ell_\eta\le 2\max\{n\in\N:\sqrt{c_n}\ge 2/\eta\}$ and
finishes the claim (b).

Finally, (c) follows easily from (b).
\end{proof}  
  \begin{lem}\label{L:2.2c}
Let $(e_i)$  be a suppression 1-unconditional,
1-subsymmetric and normalized basis.

For a normalized, basic  and weakly null sequence $(y_n)$ the following statements
 are equivalent.

\begin{enumerate}
\item[a)] $(y_n)$ has a subsequence $(z_n)$ which is strongly dominated  by $(e_i)$.

\item[b)] $(y_n)$ has a subsequence $(z_n)$ having a spreading
 model $(f_i)$ which is strongly dominated  by $(e_i)$.

\item[c)] $(y_n)$ has a subsequence $(z_n)$ for which
 there is a sequence $(\delta_n)\subset (0,\infty)$, with
 $\delta_n\searrow 0$, if $n\nearrow\infty$, so that
$$ \left\|\sum_{i=1}^\infty\alpha_i z_i\right\|\le
\max_{n\in\N} \delta_n\max_{j_1<j_2<\ldots j_n}
\left\|\sum_{i=1}^n\alpha_{j_i}e_i\right\|,
\text{ \ whenever
$(\alpha_j)\in\coo$. }$$
 \item[d)] $(y_n)$ has a subsequence $(z_n)$ for which
 there is a sequence $(\vp_n)\subset (0,\infty)$,with
 $\vp_n\searrow 0$, if $n\nearrow\infty$,
and a subsequence $(\ell_n)$ of $\N$ so that
$$ \left\|\sum_{i=1}^\infty\alpha_i z_i\right\|\le
\max_{n\in\N} \vp_n\max_{n\le j_1<j_2<\ldots j_{\ell_n}}
\left\|\sum_{i=1}^{\ell_n}\alpha_{j_i}e_i\right\|,
\text{ \ whenever $(\alpha_j)\in\coo$.}$$
\end{enumerate}
\end{lem}

\begin{proof}

\noindent (a)$\Rightarrow$(b)
Assume $(z_n)$ is a subsequence of $(y_n)$ which is strongly dominated by $(e_i)$.
 By passing to a further subsequence we can assume that $(z_n)$ has a spreading
 model $(f_i)$, being subsymmetric.
 If $\vp>0$ and
  $(\alpha_i)_{i=1}^k\subset \R$, with $\|(\alpha_i)\|_{\infty}\le \vp$, and
$\|\sum_{i=1}^k\alpha_ie_i\|\le 1$ it follows that
$$\left\|\sum_{i=1}^k\alpha_i f_i\right\|=
\lim_{n_1\to\infty}\ldots\lim_{n_k\to\infty}\left\|\sum_{i=1}^k\alpha_iy_{n_i}\right\|\le
 \Delta_{((e_i),(y_i))}(\vp),$$
which implies (b).

\noindent (b)$\Rightarrow$(a) Let $(f_i)$ be the spreading model of a subsequence $(u_i)$
 of $(y_i)$. We first note that it is enough to show that
 for a fixed $\eta>0$ we need to find a subsequence  $(v_i)$ of $(u_i)$
 and an $\vp=\vp(\eta)$ so that $\Delta_{((e_i),(y_i))}(\vp)<\eta$. Then a standard
 diagonal argument using a sequence $(\eta_n)$ decreasing to 0 will show that
 a subsequence of $(u_i)$ is  strongly dominated by $(e_i)$.

For a fixed $\eta>0$ we first choose $(\vp_n)\subset(0,1)$ so that
 $\sum_{n\in\N} \vp_n n<\eta/2$ and so that
 $\sum_{n\in\N} \Delta_{(e_i),(f_i)}(\vp_n)<\vp/4$. Since
 $(e_n)$ satisfies condition (\ref{E:2.1.1}) 
 it follows that $\ell_n=\min\{\ell\in\N:\|\sum_{i=1}^\ell e_i \|\ge
 1/\vp_{n+1}\}$ is finite for  every $n\in\N$. By passing again to a subsequence of
 we can assume that
\begin{equation}\label{L:2.2c.1}
\left\|\sum_{i=1}^{\ell_n} \alpha_i f_i\right\|\sim_2
\left\|\sum_{i=1}^{\ell_n} \alpha_i u_{k_i}\right\|,
\end{equation}
whenever $n\le k_1<k_2<\ldots k_{\ell_n}$ and $(\alpha_i)_{i=1}^{\ell_n}\subset\R$.

If $(\alpha_i)\in\coo$ with $\|\sum\alpha_i e_i\|\le 1$ and
$\max_{i\in\N}|\alpha_i|\le \vp_1$ we deduce that
\begin{align*}
\left\|\sum\alpha_i u_i\right\|&\le
 \sum_{n=1}^\infty
\left\|\sum_{i\in\N,\vp_{n+1}<|\alpha_i|\le \vp_n} \alpha_i u_i\right\|
  \le \sum_{n=1}^\infty  n\vp_n+\left\|\sum_{n\le i,\vp_{n+1}<|\alpha_i|\le \vp_n} \alpha_i u_i\right\|\\
&\le \sum_{n=1}^\infty  n\vp_n+2\left\|\sum_{n\le i,\vp_{n+1}<|\alpha_i|\le \vp_n} \alpha_i f_i\right\|
 \le \sum_{n=1}^\infty  n\vp_n+2\Delta_{((e_i)(f_i))}(\vp_n) <\eta,
\end{align*}
(for the third inequality note that since $\|\sum\alpha_i e_i\|\le 1$ it follows
 that $\#\{i:|\alpha_i|\ge\vp_{n+1}\}\le \ell_n$)
which finishes the proof of (b)$\Rightarrow$(a).

\noindent(a)$\Rightarrow$(c) Assume that $(y_i)$ is strongly dominated by
 $(e_i)$. Since we could replace $(y_n)$ by $(\tilde y_n)$ as in the remark
 after Theorem \ref{T:2.2a} and apply Proposition \ref{P:2.2b} we can assume
 that $(y_n)$ satisfies the properties (a) and (b) of
 Proposition \ref{P:2.2b}. After a renorming  we can also assume that
 $(e_i)$ dominates $(y_i)$ with constant 1, i.e. 
 that we have $\|\sum \alpha_i y_i\|\le \|\sum \alpha_i e_i\|$ for
all $(\alpha_i)\in\coo$. Choose a  strictly  to 0  decreasing sequence
 $(\vp_i)_{i=0}^\infty$ with $\vp_0=1$ and
 $\Delta_{((e_i),(y_i))}(\vp_i 2^{i+1})\le 2^{-i-1}/i$ and  define
 for $i\in\N_0$
$$K_i=\ell_{\vp_{i+1}}=
\max\left\{\#\{i:|\alpha_i|>\vp_{i+1}\}:
      (\alpha_i)\in\coo,\|\sum \alpha_i y_i\|\le 1\right\}.$$
Finally, choose for $i\in\N$
$\delta_i=2\Delta_{((e_i),(y_i))}(1)$ if $i\in\{1,2,\ldots K_0\}$ and 
 $\delta_i=1/n$ if $i\in\{K_{n-1}+1,K_{n-1}+2,\ldots, K_n\}$ for some $n\in\N$.
        
Assume now that $(\alpha_i)\in\coo$ with $\|\sum_{i=1}^\infty \alpha_i y_i \|=1$.
 Therefore
 there must be an $n\in\N_0$ for which
 $$\left\|\sum_{i\in\N,\vp_{n+1}<|\alpha_i|\le \vp_n } \alpha_i y_i \right\|\ge 2^{-n-1}
 =2^{-n-1}\left\|\sum_{i=1}^\infty \alpha_i y_i \right\|.$$
We deduce that
\begin{align*}
\left\|\sum_{i=1}^\infty \alpha_i y_i \right\|
&\le 2^{n+1}
\left\|\sum_{i\in\N,\vp_{n+1}<|\alpha_i|\le \vp_n } \alpha_i y_i \right\|\\
&=2^{n+1}\left\|\sum_{i\in\N,\vp_{n+1}<|\alpha_i|\le \vp_n } \alpha_i e_i \right\|
\left\|\sum_{i\in\N,\vp_{n+1}<|\alpha_i|\le \vp_n } \frac{\alpha_i}c y_i \right\| \\
&\left[\text{ with $c=\|\sum_{i\in\N,\vp_{n+1}<|\alpha_i|\le \vp_n } \alpha_i e_i \|$}\right]\\
&\le 2^{n+1}
 \left\|\sum_{i\in\N,\vp_{n+1}<|\alpha_i|\le \vp_n } \alpha_i e_i \right\|
 \Delta_{((e_i)(y_i))}(\vp_n 2^{n+1})\\
&\left[\text{Note that $c=\|\sum_{i\in\N,\vp_{n+1}<|\alpha_i|\le \vp_n } \alpha_i e_i \|\ge
   \|\sum_{i\in\N,\vp_{n+1}<|\alpha_i|\le \vp_n } \alpha_i y_i\|\ge 2^{-n-1}$}\right]\\
&\le\max_{j_1<j_2<\ldots j_{K_n}}\left\|\sum_{i=1}^{K_n} \alpha_{j_i} e_i\right\|
\cdot\begin{cases} 2\Delta_{((e_i),(y_i))}(1) &\text{if $n=0$}\\
             1/n &\text{if $n\ge1$}
    \end{cases}\\
&\left[\text{by choice of $K_n$}\right]\\
&=\delta_{K_n}  \max_{j_1<j_2<\ldots<j_{K_n}}\left\|\sum_{i=1}^{K_n} \alpha{j_i} e_i\right\|
 \le\max_{k\in\N} \delta_k
  \max_{j_1<j_2<\ldots< j_k}\left\|\sum_{i=1}^{k} \alpha{j_i} e_i\right\|.\\
\end{align*}
which  concludes this part of the proof.

\noindent (c)$\Rightarrow$(d)
 Let $(z_n)$ and $(\delta_n)$ be given as in (c). For
$k\in\N$, we let
 $\tilde\ell_k=\max\{\ell: \delta_\ell>2^{-k^2}\}$.

We can assume that $\tilde\ell_k\ge k$, by passing to a slower decreasing
sequence $\tilde\delta_k$, if necessary. We  also can assume that 
$(x_i)$ is a monoton basis. 

 Then it follows from (c) for
$(\alpha_i)\in\coo$ that there is an $\ell\in\N$ so that

\begin{align*}
\left\|\sum_{i=1}^\infty\alpha_iz_i\right\|
&\le\delta_{\ell}\max_{j_1<j_2<\ldots j_{\ell}}
 \left\|\sum_{i=1}^{\ell} \alpha_{j_i}e_i\right\|\\
&\le
\begin{cases}
\delta_1\max_{j_1<j_2<\ldots j_{\tilde\ell_1}}\left\|\sum_{i=1}^{\tilde\ell_1}\alpha_{j_i}e_i\right\|
&\text{if $\ell\le\tilde\ell_1$}\\
\phantom{----}\\
2^{-k^2} \max_{j_1<j_2<\ldots j_{\tilde\ell_{k+1}}}
\left\|\sum_{i=1}^{\tilde\ell_{k+1}}\alpha_{j_i}e_i\right\|
&\text{if $\ell\in[\tilde\ell_k+1,\tilde\ell_{k+1}]$}\\
\quad&\text{for some $k\in\N$.}
\end{cases}
\end{align*}
Therefore the claim follows in the case that $\ell\le \tilde\ell_1$ if we put
$\vp_1=\delta_1$ and $\ell_1=\tilde\ell_1$. 
If $\ell>\tilde\ell_1$ we can find a $k\in\N$ so that
\begin{equation*} 
\left\|\sum_{i=1}^\infty\alpha_i z_i\right\|
\le 2^{-k^2} \max_{j_1<j_2<\ldots j_{\tilde\ell_{k+1}}}
    \left\|\sum_{i=1}^{\tilde\ell_{k+1}}\alpha_{j_i}e_i\right\|\text{ and}
\left\|\sum_{i=1}^\infty\alpha_iy_i\right\|
 > 2^{-(k-1)^2} \left\|\sum_{i=1}^{\tilde\ell_{k}}\alpha_{j_i}e_i\right\|.
\end{equation*}
Thus we conclude
\begin{align*}
\left\|\sum_{i=1}^\infty\alpha_i z_i\right\|
&\le 2^{-k^2} \left[\max_{j_1<j_2<\ldots j_{\tilde\ell_{k+1}}}
    \left\|\sum_{i=k}^{\tilde\ell_{k+1}}\alpha_{j_i}e_i\right\|
  +\left\|\sum_{i=1}^{k-1}\alpha_{j_i}e_i\right\|\right] \\
&\le 2^{-k^2} \left[\max_{j_1<j_2<\ldots j_{\tilde\ell_{k+1}}}
    \left\|\sum_{i=k}^{\tilde\ell_{k+1}}\alpha_{j_i}e_i\right\|
  +\left\|\sum_{i=1}^{\tilde \ell_k}\alpha_{j_i}e_i\right\|\right] \\
&\le 2^{-k^2} \max_{k\le j_1<j_2<\ldots j_{\tilde\ell_{k+1}}}
    \left\|\sum_{i=1}^{\tilde\ell_{k+1}}\alpha_{j_i}e_i\right\|
 +2^{-k^2} 2^{(k-1)^2}\left\|\sum_{i=1}^\infty\alpha_i z_i\right\|.
\end{align*}
which implies that
\begin{align*}
\bigl[1-2^{-k^2} 2^{(k-1)^2}\bigr]
 \left\|\sum_{i=1}^\infty\alpha_i z_i\right\|
\le 2^{-k^2} \max_{k\le j_1<j_2<\ldots j_{\tilde\ell_{k+1}}}
\left\|\sum_{i=1}^{\tilde\ell_{k+1}}\alpha_{j_i}e_i\right\|.
\end{align*}
and finishes the proof of our claimed implication if
we choose $\vp_k=2^{-k^2}/[1-2^{-k^2} 2^{(k-1)^2}]$ and
$\ell_k=\tilde\ell_{k+1}$.

\noindent (d)$\Rightarrow$(a) Assume $(z_n)$, $(\vp_n)$ and $(\ell_n)$ are chosen as
required in (d).
Let $\eta>0$, choose $n_0\in\N$ so that $\vp_{n_0}\le \eta/2$ and choose
$\vp<\eta/2n_0$. 

If now $(\alpha_i)\in\coo$, $\max_i |\alpha_i|<\vp$, and
 $\|\sum\alpha_i e_i\|\le 1$ then it follows from condition
 (d) that
 \begin{equation*}
\left\|\sum_{i=1}^\infty\alpha_iz_i\right\|\le
\max_{n\in\N, n\le j_1<\ldots j_{\ell_n}} \vp_n
\left\|\sum_{i=1}^{\ell_n}\alpha_{j_i}e_i\right\|
\le \eta n_0+\max_{n\ge n_0, n\le j_1<\ldots j_{\ell_n}} \vp_n
\left\|\sum_{i=1}^{\ell_n}\alpha_{j_i}e_i\right\| <\eta, 
\end{equation*}
which implies that $\Delta_{(e_i),(z_i)}(\vp)<\eta$ and
proves the claim and finishes the proof of the Lemma.
\end{proof}

\begin{proof}[Proof of Theorem \ref{T:2.0}]

Using Lemma \ref{L:2.2c} (c)$\Rightarrow$(a)
 for the sequences $(f_i)$ and $(e_i)$ and then Lemma \ref{L:2.2c}
 (b)$\Rightarrow$(a) for the sequences  $(y_i)$ and $(e_i)$
 we can assume,
 by passing to a subsequence of $(y_n)$ if necessary, that
 $(e_i)$ strongly dominates $(y_i)$ and
 by passing to the sequence 
 $(\tilde y_n)$ as defined in the remark after Theorem \ref{T:2.2a}
we can assume that (a) and (b) of Proposition \ref{P:2.2b} are satisfied. 

Using Lemma \ref{L:2.2c} we can again pass to a subsequence, 
still denoted by $(y_n)$ for which we find sequences $(\tilde\delta_n)\subset(0,\infty)$
 and $\ell_n\subset \N$, with $\delta_n\searrow 0$ and $\ell_n\nearrow \infty$
 for $n\nearrow\infty$, so that
\begin{equation*}
\left\|\sum\alpha_i y_i\right\|\le \frac12\max_k \tilde\delta_k\max_{k\le n_1<n_2<\ldots n_k}
 \left\|\sum_{i=1}^k \alpha_{n_i} e_i\right\|,
\text{ whenever $(\alpha_i)\in\coo$.}\end{equation*}

By
 using Theorem \ref{T:2.2a} and passing to a subsequence of $(x_n)$,
if necessary, we can assume that

\begin{equation}\label{E:2.0.1}
\left\|\sum\alpha_i y_i\right\|\le \max_k\delta_k\max_{k\le n_1<n_2<\ldots n_k}
 \left\|\sum_{i=1}^k \alpha_{n_i} x_{n_i}\right\|\le 
 3 \left\|\sum\alpha_i x_i\right\|
\end{equation} 
whenever $(\alpha_i)\in\coo$.
 Thus $(x_n)$ dominates $(y_n)$ and
in order to show that the formal identity $I:x_n\mapsto y_n$
 extends to a
 strictly singular operator, let $(u_n)$ be  a seminormalized block of
$(x_n)$, write $u_n$, for $n\in\N$ as
\begin{equation*}
u_n=\sum_{i=k_{n-1}+1}^{k_n} \alpha^{(n)}_i x_i,
\text{ and let } v_n=I(u_n)\sum_{i=k_{n-1}+1}^{k_n} \alpha^{(n)}_i y_i
\end{equation*}

Using Theorem \ref{T:2.2a} and the fact that $(e_n)$ has to satisfy 
 (\ref{E:2.1.1}),
 we can assume that\\
 $\lim_{n\to\infty} \max_{k_{n-1}<i\le k_n} |\alpha^{(n)}_i|=0$,
 otherwise we could pass to an appropriate seminormalized block of
 $(u_n)$. From (\ref{E:2.0.1}) we can now easily deduce that
 $\lim_{n\to \infty} \|I(u_n)\|=0$, which proves that $I$ cannot
be an isomorphism.
\end{proof}

Theorem \ref{T:2.0} gives a necessary and sufficient condition, for the property
that a basic sequence strongly dominates an other one. Of course strong
domination is a much stronger condition then domination. Nevertheless,
 if our goal is to state a condition on the spreading models of the
 sequences $(x_n)$ and $(y_n)$ which  forces that $(x_n)$ dominates $(y_n)$ then
  strong domination of the spreading models is needed as the following remark shows.

\begin{rem} Assume that $(e_i)$ and $(f_i)$  are two normalized 1-subsymmetric 
 (1-spreading and 1-unconditional) basic
sequences, so that $(e_i)$ dominates $(f_i)$ but does not strongly dominate it.
Moreover assume that $F=[f_i:i\in\N]$ does not contain a subspace isomorphic to $c_0$.

We can therefore find  for $n\in\N$ an element $a^{(n)}=(a^{(n)}_i)\in\coo$,
so that $\max_i|a^{(n)}_i|\to 0$ if $n\to\infty$, and 
 $c',c>0$ so that
 \begin{equation}\label{ E:2.6.1}
1=\Big\|\sum_{i\in\N} a^{(n)}_if_i\Big\| \ge c \Big\|\sum_{i\in\N} a^{(n)}_ie_i\Big\|\ge c'.
\end{equation}

Now we let $(x_n)$ be the basis for the {\em Schreier space associated to }$(e_i)$, i.e.
the norm defined by
\begin{equation}\label{E:2.6.2}
\Big\|\sum a_i x_i\|=\max_{\substack{n\in\N\\ n\le i_1<\ldots i_n}}\Big\|\sum_{j=1}^n a_{i_j} e_j\Big\|
 \text{ whenever }(a_i)\in\coo.
\end{equation}

 As in the original Schreier space (where $(e_i)$ is set to be the $\ell_1$-basis) it
 is easy to see that $(e_i)$ is a spreading model of $(x_i)$.
 Let $(\tilde x_n)$ be a subsequence of $(x_n)$.
 For $n\in\N$ define $m_n=\max \supp (a^{(n)}_i)$ and $u_n=\sum_{i=1}^{m_n} a_i^{(n)}\tilde x_{m_n+i}$.
 Then $\|u_n\|\ge c'/c$ and,
 again, as in the original   Schreier space, 
 we can show that a subsequence of $u_n$ is equivalent to the $c_0$-unit basis.
 Since $F$ does not contain a copy of $c_0$ we deduce that
the map $\tilde x_n\mapsto f_n$ can not be extended to a linear bounded operator.
\end{rem}

In general the condition that  a subsymmetric and normalized  basis $(e_i)$ strongly dominates
 another basis $(f_i)$ is much stronger than the condition that $(e_i)$ dominates $(f_i)$
without $(e_i)$ being equivalent to $(f_i)$. But in  the case of $E=\ell_1$,
 we have the following.

\begin{pro}\label{P:2.7}
 Assume $(y_n)$
is a normalized basic weakly null sequence.

$(y_n)$ has a subsequence which is strongly dominated by $\ell_1$ if and only
if $(y_n)$ has a subsequence whose spreading model is not equivalent to the
unit vector  basis  of $\ell_1$.
\end{pro}

\begin{proof} By Lemma \ref{L:2.2c}  it follows immediately that $(y_n)$
is not strongly dominated by $\ell_1$ if it has a subsequence whose spreading
model is equivalent to the  unit vector basis  of $\ell_1$.

To show the converse we need to show by  Lemma\ref{L:2.2c} that if a
  1-subsymmetric   basis $(f_i)$
  is not strongly dominated by $\ell_1$ then $(f_i)$ is equivalent to
the  unit vector basis  of $\ell_1$.

Let $f^*_n$, $n\in\N$, be the coordinate functionals                                                                                      on $F$. Then $(f^*_n)$
is a  1-subsymmetric basis of the closed linear
 span of $(f^*_n)$. 

By assumption there is a $\delta_0>0$ and to for each $n\in\N$ a sequence
 $y_n=\sum_{i=1}^\infty a_i^{(n)}f_i\in[f_i:i\in\N]$ with
$0\le a_i^{(n)}<1/n$, for $i\in\N$, $\sum_{i=1}^\infty a_i^{(n)}=1$,
and $\|\sum_{i=1}^\infty a_i^{(n)}f_i\|\ge \delta_0$.

For each $n\in\N$ choose
$y^*_n=\sum_{i=1}^\infty \beta_i^{(n)} f^*_i\in [f_i^*:i\in\N]$ 
 with $0\le\beta_i^{(n)}$, for $i\in\N$,
 $\|\sum_{i=1}^\infty \beta_i^{(n)} f^*_i\|=1$, and
$y^*_n(y_n)=\sum_{i=1}^\infty a_i^{(n)}\beta_i^{(n)}\ge\delta_0/2$.

Letting $c_n=\#\{i:\beta^{(n)}_i\ge \delta_0/4\}$, $n=1,2\ldots$ it follows
from the conditions on $( a_i^{(n)})$ that
\begin{equation*}
\delta_0/2\le \sum_{i=1}^\infty  a_i^{(n)}\beta_i^{(n)}\le
c_n\frac1n+\delta_0/4,
\end{equation*}
thus $c_n\ge\frac{n\delta_0}4$. Since $(f^*n)$ is 
 1-subsymmetric it follows for all $k\in\N$ that
$\|\frac{\delta_0}4\sum_{i=1}^kf_i^*\|\le 1$ and thus
 that $f^*_i$ is equivalent to unit basis of c$_0$, from which we finally deduce
that $(f_i)$ is equivalent to the $\ell_1$-basis. 
\end{proof}
The following proposition describes another situation in which
 Theorem \ref{T:2.0} applies. Its proof can be compiled from
the techniques in  \cite{AOST}, Section 3. Nevertheless, the proof is still quite
technical, and since Theorem \ref{T:1.2} provides a generalization, we will not
give a proof here. 
Before we can state the result we need the following Definition from \cite{AOST}.

\begin{defn}\label{D:2.8}
Let $(x_i)$ be a 1-subsymmetric basic sequence. The {\em Krivine set}
of $(x_i)$ is
the set of $p$'s $(1\le p\le\infty)$ with the following property: \ For all
$\vp>0$ and $n\in\mathbb{N}$ there exists $m\in\mathbb{N}$ and
$(\lambda_k)^m_{k=1} \subset \mathbb{R}$,
such that for all $(a_i)_1^n \subseteq \mathbb{R}$,
$$\frac1{1+\vp} \|(a_i)^n_{i=1}\|_p \le \left\|\sum^n_{i=1} a_iy_i\right\| \le
(1+\vp) \|(a_i)^n_{i=1}\|_p,\text{ where}$$

$$y_i = \sum^m_{k=1} \lambda_k x_{(i-1)m+k} \quad \text{for}\quad
i=1,\ldots, n$$
and $\|\cdot\|_p$ denotes the norm of the space $\ell_p$.
\end{defn}

The proof of Krivine's theorem \cite{K} as modified by
Lemberg \cite{Le} (see also \cite{Gu}, remark II.5.14), shows
that  for every 1-subsymmetric basic sequence $(x_i)$ the Krivine set of
$(x_i)$ is non-empty. It is important to note that our definition of a Krivine
$p$ requires not merely that $\ell_p$  be block finitely representable in
$[x_i]$ but each $\ell_p^n$ unit vector basis is obtainable by means of an
identically distributed block basis.

\begin{pro}\label{P:2.9}
 Assume $X$ is a Banach space containing a normalized basic sequence $(x_n)$
 which has a
 spreading model $(e_i)$ which is not equivalent to the $\ell_1$-unit vector 
basis,
but whose Krivine set contains the number 1. 

Then there is  a   normalized basic sequence $(z_n)$ in $X$ which strongly 
dominates
 $(x_n)$.
\end{pro}

\section{The transfinite family $(\cF_\alpha)$ 
and some of its basic properties}\label{S:3}

In this section we discuss a well ordered family
$(\cF_\alpha)_{\alpha<\omega}$ of subsets of  the finite subsets of $\N$.
Its definition is similar to the definition of the Schreier family 
$(S_\alpha)_{\alpha<\omega}$  \cite{AA}.
 The Schreier set of order $\alpha$, $S_\alpha$, corresponds 
to 
our set $\cF_{\omega^{\alpha}}$, in the sense that they have the same
Cantor Bendixson index.

For every limit ordinal $\alpha<\omega_1$ we consider a sequence of sets
$(A_n(\alpha))_{n\in\N}$ so that for each $n\in\N$
 $A_n(\alpha)$ is a finite subset of $[0,\alpha)$ and
\begin{equation}\label{E:3.1}
 A_n(\alpha)\subset A_{n+1}(\alpha) \text{ for, $n\in\N$, and} 
\lim_{n\to\infty} \max A_n(\alpha)=\alpha.
\end{equation}
We call $(A_n(\alpha))_{n\in\N}$ the {\em sequence approximating
 $\alpha$}.
If for every limit ordinal $\alpha<\omega$
 (we write $\alpha\!\in\!\Lim(\omega_1)$)  $(A_n(\alpha))_{n\in\N}$ is
 a sequence approximating  $\alpha$, we call the family
$(A_n(\alpha))_{n\in\N,\alpha\in\Lim(\omega_1) }$ {\em an approximating
 family.}

Given an approximating family $(A_n(\alpha))_{n\in\N,\alpha\in\Lim(\omega_1) }$
the sets $\cF_\alpha\subset[\N]^{<\infty}$, $\alpha<\omega_1$, are defined
by  transfinite recursion as follows
 \begin{equation}\label{E:3.1a}      \cF_0=\{\emptyset\}.\end{equation}
 Assuming for some $0<\alpha<\omega_1$  the sets
$\cF_\beta\subset[\N]^{<\infty}$, with $\beta<\alpha$, are already defined we
proceed as follows.

\noindent

\begin{align}\label{E:3.2}
\cF_\alpha&=\bigl\{\{n\}\cup E: n\in\N, E\in\cF_\beta\bigr\}\cup\{\emptyset\}
 \text{ if $\alpha=\beta+1$ and}\\
 \label{E:3.3}
\cF_\alpha&=\bigl\{\ E\in[\N]^{<\infty}:  E\in \bigcup_{\beta\in A_{\min E}(\alpha)} \cF_\beta
                     \bigr\}\text{ if $\alpha\in\Lim(\omega_1)$}
\end{align}

We say that the {\em  transfinite family $(\cF_\alpha)_{\alpha<\omega}$  is defined  by the
approximating family}\\
 $(A_n(\alpha))_{n\in\N,\alpha\in\Lim(\omega_1) }$.

We first state some elementary properties of our family $(\cF_\alpha)$.

\begin{pro}\label{P:3.1}
Assume that $(\cF_\alpha)_{\alpha<\omega_1}$ is the transfinite family
associated to an approximating family 
$(A_n(\alpha))_{n\in\N,\alpha\in\Lim(\omega_1)}$.
\begin{enumerate}
\item[a)]
 For $\alpha<\omega_1$, $\cF_\alpha$ is hereditary, spreading and compact
  in $[\N]^{<\infty}$.
\item[b)]\label{P:3.2}
 For $\alpha<\omega_1$ it follows that
\begin{align*}
 \cF_{\alpha+1}&=\bigl\{\{n\}\cup E: n\in\N, n<\min E, \text{ and }
 E\in\cF_\alpha\bigr\}\cup\{\emptyset\}\\
  &=\bigl\{ E\in[\N]^{<\infty}: E\not=\emptyset, 
         E\setminus\{\min E\}\in\cF_\alpha\bigr\}\cup\{\emptyset\}.
\end{align*}
\item[c)]
 For $\alpha\le\beta<\omega_1$  there is an $m\in\N$ so that
$$\cF_\alpha\cap \bigl[ \{m,m+1,\ldots\}\bigr]^{<\infty}
\subset\cF_\beta\cap \bigl[\{m,m+1,\ldots\} \bigr]^{<\infty}.$$
\end{enumerate}
\end{pro}

\begin{proof} We can prove (a) by transfinite induction for  all $\alpha<\omega_1$ while
(b) follows from the fact that $\cF_\alpha$ is hereditary and
spreading.

To show (c)
we fix $\alpha$ and prove the claim by transfinite induction
for all $\beta>\alpha$.
Assuming the claim is true for all $\gamma<\beta$.
If $\beta=\gamma+1$ the claim follows since
 $\cF_{\gamma}\subset \cF_{\gamma+1}$.
  If $\beta$ is a limit ordinal and if $(A_n(\beta))$ is
  its approximating sequence we proceed as follows.

First we choose $n\in\N$ so that $\beta^{(n)}=\max A_n(\beta))>\alpha$, and,
 using the induction hypothesis we can find an $\ell\in\N$ so
 that
$$\cF_\alpha\cap[\{\ell,\ell+1,\ldots\}]^{<\infty}
\subset
\cF_{\beta^{(n)}}\cap[\{\ell,\ell+1\ldots\}]^{<\infty}.$$
Secondly we observe from the definition of $\cF_\beta$ it follows that
$$\cF_{\beta^{(n)}}\cap[\{n,n+1\ldots\}]^{<\infty}\subset
 \cF_{\beta}\cap[\{n,n+1\ldots\}]^{<\infty}.
$$
Therefore the claim follows by choosing $m=\max\{\ell,n\}$.
\end{proof}

In the definition of approximating families we allow the sets $A_n(\alpha)$  to have more than one element,
 contrary  to the definition of the Schreier families (see \cite{AA}), because we want to
ensure that the transfinite families are directed. 

 \begin{pro}\label{P:3.4}
Assume for every $k\in\N$
$(A^{(k)}_n(\alpha))_{n\in\N,\alpha\in\Lim(\omega_1) }$ is an approximating
 family defining the 
 transfinite family  $(\cF^{(k)}_\alpha)_{\alpha<\omega_1}$.

For each $k\in\N$, $n\in\N$ and $\alpha\in\Lim(\omega_1)$ define
 $B_n^{(k)}(\alpha)=\bigcup_{i=1}^k A^{(i)}_n(\alpha)$ and
let $\cG^{(k)}_\alpha$,  $\alpha<\omega_1$ be  defined using the approximating
 family $(B_n^{(k)}(\alpha))$.
 
 Further more define for $n\in\N$ and $\alpha\in\Lim(\omega_1)$
  $B_n(\alpha)=\bigcup_{i=1}^n A^{(i)}_n(\alpha)$ and let
  $(\cG_\alpha)_{\alpha<\omega_1}$ be the transfinite family
   defined by $(B_n(\alpha))$.

Then it follows for all $\alpha<\omega_1$
\begin{enumerate}
\item[a)] $\bigcup_{i=1}^k\cF^{(i)}_\alpha\subset\cG^{(k)}_\alpha$ for $k\in\N$.
\item[b)] For all $k\in\N$ it follows that $\cF^{(k)}_\alpha\cap
 [\{k,k+1,\ldots\}]^{<\infty}\subset\cG_\alpha$.

\end{enumerate}
\end{pro}

\begin{proof} By transfinite induction.
\end{proof}

\begin{defn}\label{D:3.4a} Let $\cA\subset[\N]^{<\infty}$ and
 $N=\{n_i:i\in\N\}\infsubset\N$, $n_i\nearrow\infty$, if $i\nearrow\infty$.
\begin{enumerate}
\item[a)]
 We call the set $\cA\cap[N]^{<\infty}$ {\em the restriction  of $\cA$ onto $N$}.
\item[b)]
We call the family 
$\cA^N=\bigl\{\{n_i:i\in E\} : E\in\cA \bigr\}$ 
 the {\em spreading of $\cA$ onto $N$.}
\end{enumerate}
\end{defn}
 Using the definition of $(\cF_\alpha)$, we
 obtain the following  recursive description
 of $\cF_\alpha^N$.

\begin{pro}\label{P:3.5}
Assume  $N\infsubset\N$, $N=(n_i)$, $n_i\nearrow\infty$
 and $\alpha<\omega_1$. Then
\begin{enumerate}
\item[a)] If $\alpha=\beta+1$
\begin{align*}
\cF_\alpha^{N}&=\bigl\{ \{n\}\cup F: n\in N\text{ and }F\in\cF_\beta^{N}\bigr\}\cup\{\emptyset\}\\
 &=\bigl\{ F\in[N]^{<\infty}\setminus\{\emptyset\}:
  F\setminus\{\min F\}\in\cF_\beta        \bigr\}\cup\{\emptyset\}.
\end{align*}

\item[b)] If $\alpha\in\Lim(\omega_1)$ and
$(A_n(\alpha))$ is the sequence approximating $\alpha$, then
 $$ \cF_\alpha^{N}=
\bigl\{ F\in[N]^{<\infty} : 
 F\in\bigcup_{\beta\in A_{\min\{i: n_i\in F\}}(\alpha)} \cF_\beta^N
\bigr\}.$$
\end{enumerate}
\end{pro}

\begin{pro}\label{P:3.6}
Assume  $M,N\infsubset\N$ and $m_0\in\N$ so that
\begin{enumerate}
\item[a)] $\#(M\cap[1,m_0])\le \#(N\cap[1,m_0])$
\item[b)] $M\cap[m_0,\infty)\subset N\cap[m_0,\infty)$
\end{enumerate}
Then it follows for $\alpha<\omega_1$ that
 $\cF_\alpha^M\cap[\{m_0,m_0+1,\ldots\}]^{<\infty}\subset\cF_\alpha^N\cap[\{m_0,m_0+1,\ldots\}]^{<\infty}$.
\end{pro}

\begin{proof} We prove the claim by transfinite induction
 on $\alpha<\omega_1$, using at each inductionstep
 Proposition \ref{P:3.5} (a) or (b).
\end{proof}

\begin{pro}\label{P:3.7}
If  $N\infsubset\N$, $N=(n_i)$, $n_i\nearrow\infty$, and
$\beta<\alpha<\omega_1$ then there is an $\ell\in\N$
so that
$$\cF_\beta^{\{n_i:i\ge \ell\}}\subset \cF_\alpha^N.$$
\end{pro}
\begin{proof} Using  Proposition \ref{P:3.1} (c) we can
choose $\ell\in\N$ so that
$\cF_\beta\cap[\{\ell,\ell+1,\ldots\}]^{<\infty}\subset
 \cF_\alpha\cap[\{\ell,\ell+1,\ldots\}]^{<\infty}$. 
Thus it follows (for the first "$\subset$'' recall that $\cF_\beta$ is spreading)
\begin{align*}
\cF_\beta^{\{n_i:i\ge \ell\}}
&=\bigl\{\{n_{i+\ell-1}:i\in E\} :E\in\cF_\beta\bigr\}
 \subset \bigl\{\{n_j:j\in E\}: E\in \cF_\beta,\text{ and } E\ge\ell\bigr\}\\  
&\subset \bigl\{\{n_j:j\in E\}: E\in \cF_\alpha, \text{ and } E\ge\ell\bigr\}
\subset\cF^N_\alpha,
\end{align*}
which finishes the proof.\end{proof}

\begin{pro}\label{P:3.8}
Assume that  $(A_n(\alpha))_{n\in\N,\alpha\in\Lim(\omega_1)}$
and $(B_n(\alpha))_{n\in\N,\alpha\in\Lim(\omega_1)}$
 are two approximating families defining the transfinite families
$(\cF_\alpha)_{\alpha<\omega_1}$ and $(\cG_\alpha)_{\alpha<\omega_1}$ 
respectively. For $\alpha<\omega_1$ and $N\infsubset \N$ there
is an $M\infsubset N$ so that
$\cG_\alpha^M\subset\cF_\alpha$.
\end{pro}

\begin{proof}
We proof the claim by transfinite induction on $\alpha<\omega_1$.
Assume that the claim is true for all $\beta<\alpha$.
If $\alpha$ is a successor it follows immediately that
the claim is true for $\alpha$.

Assume $\alpha=\sup_{\gamma<\alpha}\gamma=\sup_{n\in\N} \max B_n(\alpha)=
\sup_{n\in\N} \max A_n(\alpha)$.

Using Proposition \ref{P:3.1}(c) we can  choose an $m_i\in\N$ for each $i\in\N$ , so that
for all $\gamma\in\ B_i(\alpha)$ it follows that
 $\cF_\gamma\cap[\{m_i, m_i+1,\ldots\}]^{<\infty}\subset\cF_\alpha.$

Secondly, using the induction hypothesis, we can find $N\infsupset M_1\infsupset M_2\ldots$
 so that for all $k\in\N$ and all $\gamma\in B_k(\alpha)$ it follows that
$\cG_{\gamma}^{M_k}\subset\cF_\gamma$.  Since we can make sure that $\min M_k\ge m_k$, for $k\in\N$,
it follows that  $\cG_{\gamma}^{M_k}\subset\cF_\alpha$ for all $k\in \N$ and $\gamma\in B_k(\alpha)$.
If we finally let $M$ be a diagonal sequence of the $M_i$'s we deduce the claim
 from Proposition \ref{P:3.6}.
\end{proof}

From the property that $\cF_\alpha$ is spreading it is easy to see
that for any $L\infsubset \N$ it follows that
$\cF^L_\alpha\subset \cF_\alpha\cap[L]^{<\infty}$.
If one is willing to change the approximating family the
converse becomes true.

\begin{pro}\label{P:3.10}

Let $(\cF_\alpha)_{\alpha<\omega_1}$ be a transfinite family which is defined by
 an approximating  family $(A_n(\alpha))_{n\in\N,\alpha\in\Lim(\omega_1)}$,
and let $L\infsubset \N$.

Then there is an approximating family
  $(B_n(\alpha))_{n\in\N,\alpha\in\Lim(\omega_1)}$ defining
  a transfinite family
$(\cG_\alpha)$  for which it follows
 that for any $\alpha<\omega_1$ 
\begin{equation}\label{E:3.10.1}
\cF_\alpha\cap [L]^{<\infty}\subset\cG^L_\alpha.
\end{equation}
\end{pro}

\begin{proof}
 We write $L$ as $L=\{\ell_1,\ell_2,\ldots\ldots\}$ with
$\ell_1<\ell_2,<\ldots$.  For a limit ordinal $\alpha<\omega_1$ and an
$i\in\N$ we put
 $B_i(\alpha)= A_{\ell_i}(\alpha)$
 and show by transfinite induction that
$\cF_\alpha\cap[L]^{<\infty}\subset \cG^L_\alpha$, where
the family $(\cG_\alpha)$ is defined based on the approximating
 family $(B_i(\alpha))_{n\in\N,\alpha\in\Lim(\omega_1)}$.

 Assuming the claim to be true for all $\beta<\alpha$ the claim follows
immediately
 from Proposition \ref{P:3.5} (a)
 for $\alpha$ if $\alpha$ is a successor. If $\alpha$ is
a limit ordinal we observe that for an $F\in\cF_\alpha\cap[L]^{<\infty}$ 
 we have  $F\in\bigcup_{\beta\in A_n(\alpha)} \cF_\beta\cap[L]^{<\infty}$  with
$n=\min F$. Choosing $i\in\N_0$ so that $\ell_i= n$
 we deduce  from the induction hypothesis that
\begin{align*}
F\in\bigcup_{\beta\in A_n(\alpha)} \cF_\beta\cap[L]^{<\infty}
  \subset \bigcup_{\beta\in A_n(\alpha)} \cG_\beta^L
  = \bigcup_{\beta\in B_i(\alpha)} \cG_\beta^L
\end{align*}
which implies by Proposition \ref{P:3.5} (b) the claim.
\end{proof}

\section{The transfinite family $(\cF_\alpha)$ is universal}\label{S:4}

The main goal in this section is to prove that the family $(\cF_\alpha)$ which was introduced
in section \ref{S:3} is universal (see Theorem \ref{T:4.1} for the precise statement).

In the following Defintion we are using the transfinite family $(F_\alpha)$ to measure
the complexity of hereditary sets $\cA\subset[\N]^{<\infty}$.  As we will see
 later 
 this measure  is equivalent to the Cantor Bendixson index introduced in Section \ref{S:1}.

\begin{defn}
Consider an approximating family
 $(A_n(\alpha))_{n\in\N, \alpha\in\Lim(\omega_1)}$ defining a
 transfinite family $(\cA_{\alpha})_{\alpha<\omega_1}$.

Let $P\infsubset\N$, $\alpha<\omega_1$, and $\mathcal A\subset [\N]^{<\infty}$ be hereditary.
We say that $\mathcal A$ {\em is $\alpha$-large on $P$}, if
 for all $M\infsubset P$ there is an $N\infsubset M$ so that
 $\cF_{\alpha}^{N}\subset\cA$.

Moreover, we call
 $\I(\cA,P)=\sup\{\alpha:\cA\text{ is $\alpha$-large on $P$ }\}$
 \em{ the complexity of $\cA$ on $P$}.
\end{defn}

\noindent {\bf Remark.} Since
 the notion $\alpha$-large depends on
$(\cF_\alpha)$  which depends on the choice of
 the approximating family, we should
  have rather used   the notion $\cF_\alpha$-large instead of $\alpha$-large.
   But in Corollary \ref{C:4.3} we will show
   that the property of being $\alpha$-large is
   independent  to the underlying approximating
    family. For the results up to  Corollary \ref{C:4.3}
    we consider the approximating family and its transfinite
    family to be fixed.

\begin{pro}\label{P:4.6}  (Stabilization of $\I(\cA,P)$ with respect to $P$).

\noindent
Assume $\cA\subset[\N]^{<\infty}$ is hereditary, $P\infsubset \N$ and
 let $\alphao=\I(\cA,P)$.
Then there is a $Q\infsubset P$ so that for all $L\infsubset Q$ it follows that 
$\I(\cA,L)=\alphao$.
\end{pro}

\begin{proof} Since   $\I(\cA,P)<\alphao+1$ we deduce  that
 there is a $Q\infsubset P$ so that for all $N\infsubset Q $ it follows
that $\cF^N_{\alphao+1}\not\subset \cA$. This implies that for all $N\infsubset Q$ 
 we have $\I(\cA,N)<\alphao+1$, and thus $\I(\cA,N)\le \alphao$.
On the other hand it is clear that $\I(\cA,N)\ge I(\cA,P)$ for all $N\infsubset P$, which
implies the claim. 
\end{proof}

\begin{pro}\label{P:4.2}
Let $P\infsubset\N$, $\alpha<\omega_1$, and $\mathcal A\subset [\N]^{<\infty}$ hereditary.
\begin{enumerate}
\item[a)] If $\beta<\alpha$ and if $\cA$ is $\alpha$-large on $P$, then
 $\cA$ is   $\beta$-large on $P$.
\item[b)] If $\alpha=\beta+1$, then
\begin{align*}
\cA \text{ is $\alpha$-large on }P \iff&\forall Q\infsubset P\quad\exists
L\infsubset Q \text{ so that} \\
&(*)\begin{cases} \forall \ell\in L\quad\forall M\infsubset L\quad\exists
N\infsubset M\\
 \qquad \cF_\beta^N\subset\cA\vert_\ell:=\{E:\{\ell\}\cup
E\in\cA, \ell< E\} \\
\qquad [\text{i.e. $\cA\vert_\ell$ is $\beta$-large on
$L$}]\end{cases}
\end{align*}
\item[c)] If $\alpha\in\Lim(\omega_1)$ and $(A_n(\alpha))$ is the
approximating sequence for $\alpha$ it follows that
\begin{align*}
\cA \text{ is $\alpha$-large on }P &\iff \forall n\in\N\quad
\cA \text{ is $\max A_n(\alpha)$-large on }P\\
[&\iff(\text{by part (a) })\forall \beta<\alpha \quad \cA
\text{ is $\beta$-large on }P]
\end{align*}
\end{enumerate}

Together with (a), (c) implies that if $\I(\cA,P)<\omega_1$ it follows that the set of all $\alpha<\omega_1$ for which
 $\cA$ is $\alpha$-large  on $P$ is the closed interval 
 $[0, \I(\cA,P)]$. 

\end{pro}
\begin{proof}
(a) follows immediately from Proposition \ref{P:3.7}.

For (b) "$\Rightarrow$"  let $\cA$ be  $\alpha$-large and $Q\infsubset P$.
 Then there is an $L\infsubset Q$ with
 $\cF^L_\alpha\subset \cA$.

Since $\alpha=\beta +1$  it follows from Proposition \ref{P:3.5} that
$$\cF^L_\alpha=\bigl\{\{\ell\}\cup E:E\in\cF^L_\beta,\ell\in
L\text{ and }\ell<E\bigr\}\cup\{\emptyset\}\subset \cA.$$
 Therefore it follows for all $\ell\in L$ that
$\cF^L_\beta\subset\cA\vert_\ell$ and  therefore  it follows for any $M\infsubset L$ that 
$\cF^M_\beta\subset\cA_\ell$.

In order to prove ``$\Leftarrow$'' of (b) assume that $M\infsubset P$.
 We need to find $N\infsubset M$ so that $\cF_\alpha^N\subset \cA$.
By assumption we find an $L =(\ell_i)\subset M$ satisfying $(*)$.

By induction we can choose $n_1<n_2<\ldots $ and $N_1\supset N_2\ldots$ so that
for all $k\in\N$:
\begin{align}
\label{E:4.2.1}&\text{$n_{i+1}\in N_i$ for $i=1,\ldots k-1$}\\
\label{E:4.2.2}&\text{$\# \bigl(N_i\cap[0,n_{i+1}]\bigr)\ge i+1$ for $i=1,\ldots k-1$}\qquad\qquad\qquad\\
\label{E:4.2.3}&\text{$\cF_\beta^{N_i}\subset \cA\vert_{n_i}$, for $i=1,\ldots k$.}
\end{align}
Indeed, choose $n_1=\ell_1$ and, using the property $(*)$ (with $\ell=n_1$) 
we find
 $N_1\infsubset L $ so that $\cF_\beta^{N_1}\subset\cA\vert_{n_1}$. Then we
choose an $n_2\in N_1\cap\{n_1+1,n_1+2,\ldots\}$, large enough in
order to satisfy (\ref{E:4.2.2}), and apply $(*)$
again  (with $\ell=n_2$) to find an $N_2\infsubset N_1$  with 
$\cF_\beta^{N_2}\subset \cA\vert_{n_1}$.
 We continue in that way.

Now we claim that $\cF_\alpha^{\{n_i:i\in\N\}}\subset \cA$, which
would finish this part of the proof.
Indeed, if 
$E=\{n_{i_1},n_{i_2},\ldots n_{i_k}\}\in \cF_\alpha^{\{n_i:i\in\N\}}$,
with $n_{i_1}<n_{i_2}<\ldots<n_{i_k}$, then, by   Proposition
 \ref{P:3.5} it follows that
$\{n_{i_2},n_{i_3}\ldots n_{i_k}\}\in\cF_\beta^{\{n_i:i\in\N\}}$.                                        
Also note that by choice of $n_{i_2}$,
 $\# (N_{i_2-1}\cap[0,n_{i_2}])\ge i_2 =\# \{n_i:i\in\N\}\cap[0,n_{i_2}]$.
 Therefore, by Proposition \ref{P:3.6}
 $\{n_{i_2},n_{i_3}\ldots n_{i_k}\}\in \cF_\beta^{N_{i_2-1}}$,
 and thus, since  $\cF_\beta^{N_{i_2-1}}\subset\cF_\beta^{N_{i_1}}$, 
it follows that  
$\{n_{i_2},n_{i_3}\ldots n_{i_k}\}\in\cF_\beta^{N_{i_1}}$,
which implies by (\ref{E:4.2.3}) that 
$\{n_{i_2},n_{i_3}\ldots n_{i_k}\}\in\cA\vert_{n_{i_1}}$ and
 thus $\{n_{i_1},n_{i_2}\ldots n_{i_k}\}\in\cA$.

The claim (c)``$\Rightarrow$'' follows from (a).
In order to show  (c)``$\Leftarrow$'' let $L\infsubset P$.
 Using the assumption and part (a) we find 
$L\infsupset N_1\infsupset N_2\infsupset \ldots$ so that
 $\cF_\gamma^{N_i}\subset\cA$ for all $i\in\N$ and $\gamma\in A_i(\alpha)$. 
Then choose $N=(n_i)$ to be a diagonal sequence of
$(N_i)_{i=1}^\infty$ in such a way that for $i\in\N$ $n_i$ is in $N_i$ and 
is at least as big as the $i$-th 
   element of $N_i$.

It follows that $\cF^N_\alpha\subset\cA$. Indeed, let $E\in\cF^N_\alpha$ 
  and, thus, $E\in \cF_{\gamma_0}^N$ for some 
$\gamma_0\in A_{i_0}(\alpha)$ where $n_{i_0}=\min E$.
Note that by the choice of $N$ we have that $i_0=\#(N\cap[1,n_{i_0}])\le 
\#(N_{i_0}\cap[1,n_{i_0}])$ and 
$N\cap[n_{i_0},\infty)\subset N_{i_0}\cap[n_{i_0},\infty)$, and we
deduce from Proposition  \ref{P:3.6} that
  $E\in\cF^N_{\gamma_0}\cap[\{n_{i_0},n_{i_0}+1,\ldots\}]^{<\infty}
\subset\cF^{N_{i_0}}_{\gamma_0}\cap[\{n_{i_0},n_{i_0}+1,\ldots\}]^{<\infty}\subset \cA$,
which finishes the proof. 
\end{proof}

 Now we can conclude that  the property of being
 $\alpha$-large for a set $\cA$  does not depend on the choice
 of the approximating sequences.

\begin{cor}\label{C:4.3}
For an $\cA\subset[\N]^{<\infty}$, $P\infsubset \N$ and $\alpha<\omega_1$, the property
 of being $\alpha$-large does not depend on the choice
 of approximating family one has chosen for defining
 the sets $\cF_\beta$.
\end{cor}
\begin{proof}  Assume that $(\cG_\alpha)$ is a transfinite family defined
by another approximating family. 
We will show by transfinite induction
 that if $\cF_\alpha^N\subset\cA$ for some $N\infsubset\N$ then
there is an $M\infsubset N$ so that   $\cG_\alpha^M\subset\cA$.
Assume that our claim is true for all $\beta<\alpha$ for some
 $\alpha<\omega_1$ and let $N\infsubset\N$ be such that
$\cF_\alpha^N\subset\cA$.

If $\alpha=\beta+1$ the claim follows from the
 induction hypothesis and Proposition \ref{P:4.2} part
(b). If $\alpha=\sup_{\beta<\alpha}\beta$ the claim follows from
applying Proposition \ref{P:4.2} part (a) and (c).
\end{proof}

\begin{cor}\label{C:4.4}
 For $\cA\subset[\N]^{<\infty}$ and $P\infsubset\N$ it follows that
\begin{equation}\label{E:4.4.1}
\I(\cA,P)\le \sup_{p\in P}\sup_{L\infsubset P} \I(\cA\vert_p,L)+1.
\end{equation}
Moreover if $\I(\cA,P)$ is a limit ordinal then
 it even follows that
\begin{equation}\label{E:4.4.2}
\I(\cA,P)\le \sup_{p\in P}\sup_{L\infsubset P} \I(\cA\vert_p,L)
\end{equation}
\end{cor}

\begin{proof}
Put $\alphao=\I(\cA,P)$. If $\alphao$ is a successor our claim follows
directly from Proposition \ref{P:4.2} part (b). If $\alphao$ is a limit ordinal
 and $\beta<\alphao$ arbitrary (and thus $\beta+1<\alphao$) we conclude
  from Proposition \ref{P:4.2} part (b) that there is a $p\in P$
 and an $L\infsubset P$  so that $\I(\cA\vert_p,L)\ge \beta$.

Since $\beta<\alphao$ was arbitrary it follows that
 $\alphao\le \sup_{p\in P}\sup_{L\infsubset P} \I(\cA\vert_p,L)$.
\end{proof}

\begin{cor}\label{C:4.5} Assume that for $\cA\subset[\N]^{<\infty}$ and 
$P\infsubset\N$ we have that
\begin{equation}\label{E:4.5.1}
\forall\alpha<\omega_1\exists N_\alpha\infsubset P\quad \cF^{N_\alpha}_\alpha\subset\cA.
\end{equation}
Then there is an $L\infsubset P$ so that $[L]^{<\infty}\subset \cA$.

Therefore, if an hereditary set $\cA\subset[\N]^{<\infty}$  has the property,
 that for no $L\infsubset P$ it follows that $[L]^{<\infty}\subset\cA$,  
 the complexity of 
$\cA$ must be some countable ordinal $\alpha_0$.

We will say that the complexity of $\cA$ on $P$ is $\omega_1$  and write
 $\I(\cA,P)=\omega_1$ if (\ref{E:4.5.1}) is satisfied.
\end{cor}

\begin{proof} We first show that there is an $n\in P$, so that 
(\ref{E:4.5.1}) holds for $\cA\vert_n$ (instead of $\cA$).
Indeed, otherwise we could find for each $n\in P$ an $\alpha_n$ so that
 $\cF^{N}_{\alpha_n}\not\subset\cA\vert_n$ for  all $N\infsubset P$. Letting
$\alpha=\sup a_n$ we deduce that  $\cF^{N}_{\alpha}\not\subset\cA\vert_n$ 
  for any $n\in P$ and any  $N\infsubset P$. By Proposition \ref{P:4.2} (b) 
this would
 contradict (\ref{E:4.5.1}) for $\alpha+1$.
 Note that $n$ could have been chosen out of any given cofinite subset of $P$.

We can iterate this argument and produce a strictly increasing  sequence 
$(n_i)\subset P$
 so that for every $k\in\N$  (\ref{E:4.5.1}) holds for 
$\cA\vert_{n_1,n_2,,\ldots n_k}=
 \{A\subset \N: \{n_1,\ldots n_k\}\cup A\in\cA\}$ holds. This implies that
 that $[\{n_i:i\in\N\}]^{<\infty}\subset \cA$.
\end{proof}

\begin{cor}\label{C:4.5a}
 For $\alpha<\omega_1$ and $L\infsubset \N$ it follows that
 $\I(\cF_\alpha, L)=\alpha$.
\end{cor}
\begin{proof} Since $\cF_\alpha^L\subset \cF_\alpha\cap[L]^{<\infty}$ it is 
clear
 that $\I(\cF_\alpha, L)\ge\alpha$.
Assume that for some $\beta>\alpha$ and some 
$N=\{n_1^{(1)},n_2^{(1)}\ldots\}\infsubset \N$ it
follows  that 
$\cF_{\beta}^N\subset \cF_{\alpha}$.

By Propostion \ref{P:4.2}  (a) we can assume that
 $\beta=\alpha+1$ and we claim that it would
 follow that there is a family $(N_\beta)_{\alpha<\beta<\omega_1}$
 of infinite subsets of $\N$,  with $N_\gamma\setminus N_\beta$ being finite, if $\gamma<\beta$, so that
\begin{equation}\label{E:4.5a.2}
\cF_{\beta}^{N_\beta}\subset \cF_{\alpha}.
\end{equation}
Using Corollary \ref{C:4.5} this would imply that for some $L\infsubset\N$ 
so that
 $[L]^{<\omega}\subset \cF_\alpha$, contradicting the compactness
 of $\cF_\alpha$.

We will show the existence of $N_\beta$  by transfinite induction of $\beta>\alpha$.
For $\beta=\alpha+1$ $N_\beta$ exists by assumption.
If $\beta=\gamma+1$ and if $N_\gamma=\{n^{(\gamma)}_1,n^{(\gamma)}_2,\ldots\}\infsubset\N$ is
as in (\ref{E:4.5a.2}) we choose
 $N_{\gamma+1}=\{n^{(\alpha+1)}_{n^{(\gamma)}_i}: i\in\N\}$ and observe that
 (note that $\cF^{N_\gamma}_{\gamma+1}\subset \cF_{\alpha+1}$ since
 $\cF^{N_\gamma}_{\gamma}\subset \cF_{\alpha}$)
\begin{align*}
\cF^{N_{\gamma+1}}_{\gamma+1}&=
 \Bigl\{\{n^{(\alpha+1)}_{n^{(\gamma)}_i}:i\in E\}: E\in \cF_{\gamma+1}\Bigr\}
=\Bigl\{\{n^{(\alpha+1)}_m:m\in F\} : F\in \cF^{N_\gamma}_{\gamma+1}\Bigr\}\\
 &\subset \Bigl\{\{n^{(\alpha+1)}_m:m\in F\} : F\in \cF_{\alpha+1}\Bigr\}
=\cF^{N_{\alpha+1}}_{\alpha+1}\subset \cF_\alpha.
\end{align*}
If $\beta\in\Lim(\omega_1)$ and if $(A_i(\beta))$ is  the sequence approximating
 $\beta$ we first note that we can choose infinite subsets
$\N\infsupset M_1\infsupset M_2\ldots$ so that for all $i\in\N$ it follows that
 $\bigcup_{\gamma\in A_i(\beta)}\cF^{N_i}_\gamma\subset \cF_\alpha$.
Then we choose $N_\beta=\{n^{(\beta)}_i: i\in\N\}$ with $n_i^{(\beta)}$ being
 the $i$-th element of $M_i$, for $i=1,2\ldots$. We deduce then the claim 
in this case from Proposition \ref{P:3.6}.
\end{proof}
Using  our results on the family $(\cF_\alpha)_{\alpha<\omega}$ we can now show 
the relation ship between $\I$ and the Cantor Bendixson index.

\begin{lem}\label{L:5.5} Assume that $\cA\subset [\N]^{<\infty}$ is hereditary and compact.
\begin{enumerate}
\item[a)] Define the {\em expansion of $\cA$} by 
$\Ep(\cA)=\bigl\{\{n\}\cup A: A\in\cA\bigr\}\cup\{\emptyset\}$.
(Note that $\cF_{\alpha+1}=\Ep(\cF_\alpha)$ for $\alpha<\omega_1$ and that
$F\in\Ep(\cA)$ if and only if $F=\emptyset$ or $F\setminus\{\min F\}\in\cA$.)\\
 For $\alpha\le\CB(\cA)$
 it follows that
 $\Ep\bigl(\cA^{(\alpha)}\bigr)=\bigl(\Ep(\cA) \bigr)^{(\alpha)}$.
\item[b)] Assume that there is a sequence $(\ell_n)\subset\N$ with $\lim_{n\to\infty}\ell_n=\infty$ and
 for $n\in\N$ hereditary and compact sets 
 $\cA_n\subset[\{\ell_n,\ell_n+1,\ell_n+2,\ldots\}]^{<\infty}$
 so that $\alpha_n=\CB(\cA_n)$  strictly increases to some $\alpha<\omega_1$, and assume that
  $\cA=\bigcup\cA_n$. Then it follows that $\CB(\cA)=\alpha+1$.
\item[c)] For any $n\in\N$ and any $\alpha<\omega_1$ it follows that
 $\bigl(\cA|_n\bigr)^{(\alpha)}=\bigl(\cA^{(\alpha)}\bigr)|_n$
\end{enumerate}
\end{lem}
\begin{proof} (a) First assume that $\alpha=1$ and, thus, that $\CB(\cA)\ge 1$. Clearly, 
 $\Ep\bigl(\cA^{(1)}\bigr)$ as well as $\bigl(\Ep(\cA) \bigr)^{(1)}$ contain $\emptyset$ as an element.
For $F\in[\N]^{<\infty}$, $F\not=\emptyset$, we observe that
\begin{align*}
F\in\bigl(\Ep(\cA) \bigr)^{(1)}
&\iff\exists N\infsubset\N, \min N>\max F\,\forall n\in N \quad F\cup\{n\}\in \Ep(\cA)\\
&\iff\exists N\infsubset\N, \min N>\max F\,\forall n\in N \quad \bigl(F\cup\{n\}\bigr)\setminus \{\min F\cup\{n\}\}\in\cA\\
 &\iff\exists N\infsubset\N, \min N>\max F\,\forall n\in N \quad \bigl(F\setminus\{\min F\}\bigr)\cup\{n\}\in\cA\\
&\iff  F\setminus\{\min F\}\in \cA^{(1)}\iff F\in \Ep\bigl(\cA^{(1)}\bigr).
\end{align*}
For general $\alpha\le\CB(\cA)$ the claim now follows easily by transfinite induction.

\noindent (b)  Note that for  $\beta<\alpha$ we can choose $m\in\N$ so that $\alpha_m>\beta$ and thus
$\emptyset\in\cA_m^{(\beta)}\subset \cA^{(\beta)}$.
Since $\alpha$ is a limit ordinal we deduce that
$\emptyset\in  \cA^{(\alpha)}$ and, thus,  that $\CB(\cA)\ge \alpha+1$.

On the other hand note that for any $m\in\N$  it follows that
\begin{equation*}
\cA^{(\alpha_m)}=\Bigr(\bigcup_{i=1}^m\cA_i\Bigr)^{(\alpha_m)}\cup \Bigr(\bigcup_{i=m+1}^\infty\cA_i\Bigr)^{(\alpha_m) }
 \subset [\{\ell_m,\ell_m+1,\ldots\}]^{<\infty}
\end{equation*}
and, thus, that $\cA^{(\alpha)}\subset\bigcap_{m\in\N}[\{\ell_m,\ell_m+1,\ldots\}]^{<\infty}=\{\emptyset\}$, which
implies that $\CB(\cA)\le \alpha+1$.

\noindent  To prove (c)  let first $\alpha=1$ and $n\in\N$. 
For $F\in[\N]^{<\infty}\setminus\{\emptyset\}$, $\min F>n$,
it follows that
\begin{align*}
F\in \bigl(\cA|_n\bigr)^{(1)}&\iff \exists M\infsubset\N, \min M>n\forall m\in M\quad
 \{n\}\cup F \cup   \{m\}\in \cA\\ 
&\iff   \{n\}\cup F\in\cA^{(1)}\iff F\in \bigl(\cA^{(1)}\bigr)|_n
\end{align*}
For general $\alpha$ we  conclude the claim by transfinite induction.
\end{proof}
\begin{cor}\label{C:5.6} For $\alpha<\omega_1$, a hereditary and compact $\cA\subset [\N]^{<\infty}$,
and an $P\infsubset \N$  it follows that
\begin{align}\label{E:6.6.1}
 \cA \text{ is $\alpha$-large on }P
\iff \forall Q\infsubset P \quad\CB(\cA\cap[Q]^{<\infty})\ge \alpha+1.
\end{align} 
\end{cor}

\begin{proof} In order to show ``$\Rightarrow$'' it is enough to observe that
 $\CB(\cF^N_\alpha)=\alpha+1$  for any $N\infsubset$ which follows by transfinite induction
 on $\alpha$ using (a) (in the successor case) and (b) (in the case of limit ordinals) of Lemma
 \ref{L:5.5}.

We also show ``$\Leftarrow$'' by transfinite induction on $\alpha<\omega_1$ and assume that
 the implication is true for all $\tilde\alpha<\alpha$.

Assume that $\cA\subset[\N]^{<\infty}$ is compact and hereditary  
 so that for all $Q\infsubset P$ it follows that $\CB(\cA\cap[Q]^{<\infty})\ge \alpha+1$.

If $\alpha=\beta+1$, it is 
 by the induction hypothesis and  by Proposition \ref{P:4.2}  enough to show:

\noindent{\bf Claim.}
$\forall Q\infsubset P\,\exists L\infsubset Q\,\forall \ell\in L \,\forall K\infsubset L,
 \ell<\min K\quad
\CB(\cA\cap[K]^{<\infty})\ge \beta+1$.

Assume the claim is not true and choose a $Q\infsubset P$
 so that for all $L\infsubset Q$ there is an $\ell\in L$ and a $K\infsubset L$ so that
 $\CB(\cA\cap[K]^{<\infty})\le \beta$.

We  first put  $L_1=Q$  and then choose $\ell_1\in L_1$  and $K_1\infsubset L_1$, with $\min K_1>\ell_1$
 so that $\CB(\cA|_{\ell_1}\cap[K_1]^{<\infty})\le \beta$. Then we let $L_2=K_1$
 and choose an $\ell_2\in L_2$ and a $K_1\infsubset L_2$ with $\min K_2>\ell_2$ so that
 $\CB(\cA|_{\ell_2}\cap[K_2]^{<\infty})\le \beta$. 
We can continue in this way and eventually get a strictly increasing sequence $L=(\ell_i)$ and
 $(L_i)$ with $Q=L_1\infsupset L_2\infsupset\ldots $ so that
 $\CB(\cA|_{\ell_i}\cap[L_{i+1}]^{<\infty})\le \beta$.
  Thus it follows  for each $i\in\N$ that
\begin{equation*} 
\bigl(\cA|_{\ell_i}\cap [L]^{<\infty}\bigr)^{(\beta)}=
 \bigl(\cA|_{\ell_i}\cap [\{\ell_{i+1},\ell_{i+1}+1,\ldots\}]^{<\infty}\bigr)^{(\beta)}
\subset \bigl(\cA|_{\ell_i}\cap [L_{i+1}]^{<\infty}\bigr)^{(\beta)}=\emptyset,
\end{equation*}
which implies by Lemma \ref{L:5.5} (c) that $\cA^{(\beta)}\cap [L]^{<\infty}$ must be 
  finite and, thus, that \\
$\CB(\cA\cap [L]^{<\infty})\le\beta+1$, contradicting the assumption. 
This proves the claim and the induction step in the case that $\alpha$ is a successor.

If $\alpha$ is a limit ordinal and $\CB(\cA\cap[Q]^{<\infty})\ge \alpha+1$ for all
 $Q\infsubset P$ it follows that for all $\beta<\alpha$ we have $\CB(\cA\cap[Q]^{<\infty})\ge\beta$,
and, thus, by our induction hypothesis that $\cA$ is $\beta$-large on $P$, for all
 $\beta<\alpha$, which implies, by Proposition  \ref{P:4.2} (c), that 
$\cA$ is $\alpha$-large on $P$.
\end{proof}

Using the equivalence of the strong Cantor Bendixson index and the concept
 of $\alpha$-largeness and using Corollary \ref{C:4.5} we can prove Proposition \ref{P:1.1} of Section \ref{S:1}.

\begin{proof}[Proof of Proposition \ref{P:1.1}.]

Assume that $X$ is a Banach space with a seminormalized basis $(x_i)$
  and let $(z_n)$ 
 another seminormalized basic sequence with $\gamma_0((z_n),(x_{k_n}))=\omega_1$.

By Corollary \ref{C:5.6} we deduce for $1\le c<\infty$ and $\gamma<\omega_1$ that
\begin{align*}
c((z_n),(x_n),\gamma)=
\inf\Bigl\{ c\ge 1: \exists N\infsubset\N\quad\cC((z_n),(x_n),c)\text{ is $\gamma$-large on N}\Bigr\}<\infty.
\end{align*}
Since $c((z_n),(x_n),\gamma)$ is non decreasing in $\gamma$ we deduce from the uncountability of
$[0,\omega_1)$ that there is a $c_0>0$ so that
 $c((z_n),(x_n),\gamma)\ge c_0$ for all $\gamma<\omega_1$. 
 But this means that for any $\gamma<\omega_1$ the set
  $\cC((z_n),(x_n),c_0/2)$ is $\gamma$ large on some set $N_\gamma\infsubset\N$.
 From Corollary \ref{C:4.5} we deduce therefore that there is an  $L=(\ell_i)\infsubset\N$
 so that $[L]^{<\infty}\subset\cC((z_n),(x_n),c_0/2)$. By passing to a subsequence
 of $(z_n)$ we might simply assume that $L=\N$.

 Therefore we can choose for any $n\in\N$ a sequence $m^{(n)}_1< m^{(n)}_2<\ldots m^{(n)}_n$ in $N$
so that $(z_i)_{i=1}^n$ is $c_0/2$ equivalent to $(x_{m^{(n)}_i})_{i=1}^n$. Passing possibly to
  a  subsequence of $(z_n)$ and  having possibly to redefine
 the $m^{(n)}_i$'s (note for a fixed sequence $(n_k)\subset\N$ we could change the 
 choice of $m^{(n)}_i$ in such a way that $m_i^{(n)}=m_i^{(n_k)}$ if $i\le n$ and $n\in[n_{k-1}+1,n_k]$)
 we can assume one of the following two cases happens.

\noindent Case 1. There is a sequence $(m_n)$ so that $m^{(n)}_i=m_i$ for all $n\in\N$ and $i\le n$.

In this case it follows that $(z_i)$ is isomorphically equivalent to $(x_{m_i})$.

\noindent Case 2. For any $i\in\N$ it follows that $\lim_{n\to\infty} m^{(n)}_i=\infty$.

In the second case it follows that $(z_i)$ is equivalent to a spreading model of a subsequence of
 $(x_n)$. This proves the first part of  Proposition \ref{P:1.1}. 

In order to deduce the second part 
 we assume that for all block bases $(z_n)$  of $(x_n)$ it follows  
 $\gamma_0((z_n),(x_n))=\gamma_0((x_n),(z_n))=\omega_1$  and fix a block basis
 $(z_n)$. By the first part of the proof a subsequence of $(z_n)$ could be equivalent to
 a subsequence of $(x_n)$, then  we are done.
 Otherwise $(z_n)$ is equivalent to a spreading model of a subsequence of $(x_n)$ which
means  in particular that $(z_n)$ is subsymmetric. Now we change the roles of $(x_n)$ and $(z_n)$,
 use the assumption $\gamma_0((x_n),(z_n))=\omega_1$, and go again through the arguments
 of the first part of the proof and  observe that 
since  $(z_i)$ is subsymmetric both cases collaps to one and that
  a subsequence of $(x_n)$ is isomorphically equivalent to a subsequence of $(z_i)$.
This proves that $X$ is  a space of Class 1.
 The other direction of the stated equivalence is trivial.
\end{proof}

We are now in the position to state and prove Theorem \ref{T:4.1} concerning
the universality of the transfinite families.

\begin{thm}\label{T:4.1}
Let $\cA\subset [\N]^{<\infty}$ be  not empty and  hereditary,
 and assume that $\alphao=\I(\cA,P)<\omega_1$. 
Then there  is an approximating family
$(B_n(\alpha))_{n\in\N,\alpha\in\Lim(\omega_1)}$ defining   the  
transfinite family
 $(\cG_\alpha)_{\alpha<\omega_1}$,  and 
there is an $L\infsubset P$ so that
\begin{equation}\label{E:4.1.1}
        \cG_{\alphao}^L\subset\cA\cap[L]^{<\infty}\subset\cG_\alphao\cap[L]^{<\infty}.
\end{equation}
\end{thm}

\begin{proof} 
Let $(\cF_\alpha)$ be a transfinite family being chosen a priori.
We will
 prove the claim by transfinite induction
 for all $\alphao<\omega_1$.

 If $\alphao=0$, we deduce that $L=\{\ell\in P: \{\ell\}\not\in\cA\}$ is infinite
and, thus, since $\cA$  is hereditary and not empty it follows that
\begin{equation*}
\cA\cap [L]^{<\infty}=\{\emptyset\}=\cF_0.
\end{equation*} 
Assume the claim to be true
 for all hereditary $\tilde\cA\subset[\N]^{<\infty}$ with
  $\I(\tilde\cA,P)<\alphao$, where $\alphao\ge1$.  Let
$\cA\subset[\N]^{<\infty}$ be hereditary with $\I(\cA,P)=\alphao$.
 By passing to a
subsequence  of $P$, if necessary,  we can assume that  $\I(\cA,L)=\alphao$ for
all $L\infsubset P$ (we are using Proposition \ref{P:4.6}). 
 Since $\cA$ is not $\alphao+1$-large on $P$ we deduce from
 Proposition \ref{P:4.2} part (b) that there is a $Q\infsubset P$ so
that
\begin{equation}\label{E:4.7.2}
\forall M\infsubset Q \text{ }\exists m\in M\quad \cA\vert_m\text{ is not
$\alphao$-large on $M$}.
\end{equation}
We start by applying (\ref{E:4.7.2}) to $M_1=Q$ and
find an $m_1\in M_1$ for which
 $\beta_{m_1}=\I(\cA\vert_{m_1},M_1)<\alphao$
 (recall that by Proposition \ref{P:4.2} the set  of ordinals 
  $\alpha$ for which $\cA\vert_{m_1}$ is $\alpha$-large
 is a closed interval).
 By the induction hypothesis we can
 find an approximating family
 $(B_n^{(1)}(\gamma))_{n\in\N,\gamma\in\Lim(\omega_1)}$
 which defines
 a transfinite family    $(\cG^{(1)}_\gamma)_{\gamma<\omega_1}$
 and an $M_2\subset M_1$
so that
$\cA\vert_{m_1}\cap[M_2]^{<\infty}\subset
\cG^{(1)}_{\beta_{m_1}}\cap[M_2]^{<\infty}$.

Since $\cA\vert_{m_1}$ is $\beta_{m_1}$-large on $M_1$
(which does not depend on the choice of the approximating
family)
 we also
 can require that
$\bigl(\cG^{(1)}_{\beta_{m_1}}\bigr)^{M_2}\subset \cA\vert_{m_1}$.

By repeating this argument we find
an increasing sequence $(m_i)_{i\in\N}$, 
sets $M_i\infsubset\N$, for $i\in\N$, a sequence of ordinals
 $(\beta_{m_i})_{i\in\N}\subset[0,\alphao) $, and approximating families
 $(B_n^{(i)}(\gamma))_{n\in\N,\gamma\in\Lim(\omega_1)}$
defining transfinite families
  $(\cG^{(i)})_{\gamma<\omega_1}$, for $i\in\N$, 
so that for all $i\in\N$

\begin{equation}\label{E:4.7.3}
m_i\in M_i, M_{i+1}\infsubset M_i,\text{ and }m_i<\min  M_{i+1},
\end{equation}
\begin{equation}\label{E:4.7.4}
\bigl(\cG^{(i)}_{\beta_{m_i}}\bigr)^{M_{i+1}}
\subset\cA\vert_{m_i}\cap[M_{i+1}]^{<\infty}\subset
\cG^{(i)}_{\beta_{m_i}}\cap[M_{i+1}]^{<\infty},
\end{equation}
Putting $M=\{m_1,m_2,\ldots\}$ we deduce from (\ref{E:4.7.3}) and
(\ref{E:4.7.4}) that for all $i\in\N$ 
\begin{equation}\label{E:4.7.6}
\bigl(\cG^{(i)}_{\beta_{m_i}}\bigr)^{\{m_{i+1},m_{i+2},\ldots\}}
\subset\cA\vert_{m_i}\cap[M]^{<\infty}\subset
\cG^{(i)}_{\beta_{m_i}}\cap[M]^{<\infty}.
\end{equation}
(for the second ``$\subset$'' recall that we defined $\cA\vert_m$ in such a way that
 $\cA\vert_m\subset[\{m+1,m+2,\ldots\}]^{<\infty}$)
which  implies that for any $m\in M$ and any
 $\tilde M\infsubset M$ we have
$\I(\cA\vert_m,\tilde M)=\beta_m$ and, thus, by
 Corollary \ref{C:4.4}, it follows
that
\begin{equation}\label{E:4.7.7}
 \sup_{m\in M,m \ge k} \beta_m+1=\alphao,\quad \text{ for all }k\in\N.
\end{equation}

To finish the proof we distinguish between
 the  case that $\alphao$ is a successor and
the case that  $\alphao$ is a limit ordinal.

If $\alphao=\gamma+1$ we deduce from (\ref{E:4.7.7}) that the set
$\tilde L=\{m\in M:\beta_m=\gamma\}$ is infinite. Using
Proposition \ref{P:3.4} (b) we can find an
approximating family $(B_n(\alpha))_{n\in\N,\alpha<\omega_1} $
 so that for any $\alpha<\omega_1$ and any $i\in\N$ it
follows that  $\cG^{(i)}_\alpha\cap[\{i,i+1\ldots\}]^{<\infty}\subset\cG_\alpha$, 
 where $(\cG_\alpha)$ is defined
 by   $(B_n(\alpha))_{n\in\N,\alpha<\omega_1} $.
We therefore deduce that
\begin{align*}
\cA\cap[\tilde L]^{<\infty}&=
\bigl\{\{\ell\}\cup E:\ell\in\tilde L, \ell< E
\text{ and }E\in\cA\vert_\ell\bigr\}\cap[\tilde L]^{<\infty}\cup\{\emptyset\}\\
 &\subset\!
\bigl\{\{\ell\}\cup E:\exists i\in\N\quad \ell=m_i\in\tilde L, \ell< E,
 E\in\cG_\gamma^{(i)} \bigr\}
\!\cap\![\tilde L]^{<\infty}\!\cup\!\{\emptyset\} \text{\ (by (\ref{E:4.7.6}))}\\
&\subset \!\bigl\{\{\ell\}\cup E: \ell\in\tilde L,\ell< E
\text{ and }E \!\in\!\cG_\gamma\bigr\}\cap[\tilde L]^{<\infty}\cup\{\emptyset\}
          \!=\!\cG_\alphao\!\cap\![\tilde L]^{<\infty}      \text{\ (since $\!i \!\le \! m_i$).}
\end{align*}

If $\alphao\in\Lim(\omega_1)$ we also define
 the approximating  family
$(B_n(\alpha))_{n\in\N,\alpha<\omega_1} $
as in Proposition \ref{P:3.4} (b), but  add in the case
 of $\alpha=\alphao$ the ordinals $\beta_{m_1}+1,\beta_{m_2}+1,\ldots \beta_{m_n}+1$ to
 the set  $B_n(\alphao)$ (still denoting it
  $B_n(\alphao)$).
 The transfinite family defined by  $(B_n(\alpha))_{n\in\N,\alpha<\omega_1} $
 is denoted by $(\cG_\alpha)$.
 We put $\tilde L=M=\{m_1,m_2,\ldots\}$.

 Now if $E\in\cA\cap[\tilde L]^{<\infty}$, $E\not=\emptyset$,
 we write $E=\{m\}\cup F$ with
 $m=m_n=\min E$ and $F\in\cA\vert_{m_n}\subset\cG^{(n)}_{\beta_{m_n}}$ 
 (by (\ref{E:4.7.6})).
 From the definition   of the family $(\cG_\alpha)_{\alpha<\omega_1}$ we conclude
 that $F\in\cG_{\beta_{m_n}}$
 and thus $E\in\cG_{\beta_{m_n}+1}$. Since
 $n\le m_n\le E$ and  since
  $\beta_{m_n}+1\in B_n(\alphao)$  we deduce
that $E\in\cG_\alphao\cap[\tilde L]^{<\infty}$.

Therefore we derive in both cases ($\alphao$ being a successor and $\alphao$ being a 
 limitordinal) that
\begin{equation}\label{E:4.7.8}
\cA\cap[\tilde L]^{<\infty}\subset\cG_\alphao\cap[\tilde L]^{<\infty}.
\end{equation}

On the other hand it follows from the definition of $\alphao$ and from
Proposition \ref{P:4.2} that $\cA$ is $\alphao$-large (which by Corollary \ref{C:4.3}
 does not depend on the transfinite family).   
We can therefore chose an $L\infsubset \tilde L$ so that
 $\cG_{\alphao}^L\subset\cA $ which implies that  
\begin{equation*}
\cG_{\alphao}^L\subset
\cA\cap[L]^{<\infty}\subset\cG_\alphao\cap[L]^{<\infty}
\end{equation*}
and finishes the proof of Theorem \ref{T:4.1}.
\end{proof}

We will have to apply Theorem \ref{T:4.1} not only for one 
$\cA\subset[\N]^{<\infty}$ but for a sequence $(\cA_n)$ simultaneously.
 Therefore we  need the following reformulation.

\begin{cor}\label{C:4.10}
Assume we are given a $P\infsubset N$, a sequence $(\cA_{\ell})_{\ell\in\N}$ of nonempty
 and hereditary subsets of $[\N]^{<\infty}$, an increasing sequence $(\ell_k)\subset\N$,
and a sequence of ordinals $(\alpha_\ell)$ so that 
\begin{equation}\label{E:4.10.1}
\alpha_\ell\ge \I(\cA_\ell, Q),\text{ whenever $\ell\in\N$ and }Q\infsubset P.
\end{equation}
Then there is a transfinite family $(\cG_\alpha)_{\alpha<\omega_1}$ and $K\infsubset P$,
$K=\{k_1,k_2,\ldots\}$, $(k_i)$ strictly increasing, so that for all $n\in\N$ and $\ell\le\ell_n$
\begin{equation}\label{E:4.10.2}
\cA_\ell\cap [\{k_n,k_{n+1},\ldots\}]^{<\infty}\subset G^K_{\alpha_\ell}.
\end{equation}
\end{cor}
\begin{proof} 
 We first use Proposition \ref{P:4.6}
 and an easy diagonalization argument to assume that $\I(\cA_\ell, Q)=\I(\cA_\ell,P)$ for all
  $Q\infsubset P$ (note that (\ref{E:4.10.1}) stays valid if we pass to subsequences).

For  $i\in\N$  we then apply Theorem \ref{T:4.1}  to  each of
 $\cA_\ell$, ${\ell\le \ell_i}$, in order to get
an approximating family 
$(B^{(i)}_n(\gamma))_{n\in\N,\gamma\in\Lim(\omega_1)}$ with associated 
 transfinite families $(\cG^{(i)}_\alpha)_{\alpha<\omega_1}$ and an $L_i\infsubset P$
 so that for any $\ell\le \ell_i$ it follows that
 $L_i\infsubset L_{i-1}$ (with $L_0=P$) and
  $\cA_\ell\cap[L_i]^{<\infty}\subset \cG_{\alpha_\ell}$.

Now define for $n\in\N$ and $\gamma\in\Lim(\omega_1)$ as in Proposition 
\ref{P:3.4}, i.e. $B_n(\gamma)=\bigcup_{j=1}^n B_n^{(j)}(\gamma)$, let
 $(\cH_\alpha)$ be the transfinite family associated to $B_n(\gamma)$,
and Let $K=\{k_1,k_2\ldots\}$ be a diagonal sequence of
the $L_i$'s. Then we deduce
that for any $m\in\N$ and any $\ell\le \ell_m$ it follows that
\begin{align*}
\cA_\ell\cap[\{k_{m},k_{m+1}\ldots\}]^{<\infty}
&=\cA_\ell\cap [L_m]^{<\infty}\cap [\{k_{m},k_{m+1},\ldots\}]^{<\infty}\\
&\subset \cG_{\alpha_\ell}^\ell\cap[\{k_{m},k_{m+1},\ldots\}]^{<\infty}
 \subset\cH_{\alpha_\ell}\cap [K]^{\infty}.
\end{align*}
By  Proposition \ref{P:3.10} we then choose $(G_\alpha)_{\alpha<\omega_1}$  so that 
 $\cH_\alpha\cap[K]^{<\infty}\cG^K_\alpha$, for $\alpha<\omega_1$.
\end{proof}

\section{Some Consequences of Theorem \ref{T:4.1}}\label{S:5}

Using Theorem \ref{T:4.1} we deduce the following generalization
of Ramsey's theorem for finite sets.

\begin{cor}\label{C:5.1}
Let $\cF\subset[\N]^{<\infty}$ be hereditary and let $P\infsubset \N$.

If $\cF=\cA\cup\cB$ with $\cA$ and $\cB$ also being hereditary then
there is a $Q\infsubset P$ so that

\begin{equation}\label{E:5.1.1}
\max(\I(\cA,Q),\I(\cB,Q))=\I(\cF,Q)
 \end{equation}
\end{cor}

\begin{proof} By passing to an infinite subsequence of $P$ and using Proposition \ref{P:4.6}, if necessary, 
we can assume that there are ordinals  $\alpha_0,\beta_0$ and $\gamma_0$ so that
 for all $\tilde P\infsubset P$
\begin{equation*}
\alpha_0=\I(\cA,\tilde P),\quad\beta_0=\I(\cB,\tilde P),\text{ and }\gamma_0=\I(\cF,\tilde P).
\end{equation*}
We need to show that $\max(\alpha_0,\beta_0)=\gamma_0$.

Assume that this is not true and, thus,  assume that $\max(\alpha_0,\beta_0)<\gamma_0$.

 By Theorem \ref{T:4.1} 
 and the fact that $\I(\cA,P)$ and $\I(\cB,P)$ are stabilized in the sense of
 Proposition \ref{P:4.6} 
 we find a transfinite family
 $(\cG_{\alpha})_{\alpha<\omega_1}$ and
an $L\infsubset P$ so that $\cA\cap[L]^{<\infty}\subset \cG_{\alpha_0}$
  and $\cB\cap[\tilde L]^{<\infty}\subset \cG_{\beta_0}$ and therefore, it would follow from
 Corollary \ref{C:4.5a}
 that $\I(\cA,L)\le \I(\cG_{\max(\alpha_0,\beta_0)},L)= \max{(\alpha_0,\beta_0)}<\gamma_0$ 
which is a contradiction.
\end{proof}

 We introduce the following "addition" of subsets of $[\N]^{<\infty}$
\begin{defn}\label{D:5.7} For $\cA,\cB\subset[\N]^{<}$ we define
 $\cA\sqcup \cB =\{ A\cup B: A\in\cA,B\in\cB\}$
\end{defn}

\begin{rem} At first sight one might believe 
that  $\I(\cF_\alpha\sqcup\cF_\beta,N)=\alpha+\beta$. But this cannot be true 
 since on one hand addition on the ordinal numbers is not commutative, on the other 
hand it is clear that $\cA\sqcup\cB=\cB\sqcup\cA$, for any $\cA,\cB\subset[\N]^{<\infty}$.
 For example it is  easy to see that $\cF_{\omega+1}\sqcup\cF_\omega=2\omega+1\not=2\omega=\omega+1+\omega$.

One could define the following "commutative addition of ordinal numbers":
\begin{equation*} 
\alpha\sqcup\beta=\I(\cF_\alpha\sqcup\cF_\beta,\N).
\end{equation*}
It might be interesting to  determine the properties of
 this  binary operation and compare it with the addition of ordinal numbers.
Nevertheless it is easy to prove by transfinite induction on $\beta$ that for all
 $\alpha,\beta<\omega_1$ it follows
that 
\begin{equation}\label{E:5.7.1}
\cF_{\alpha+\beta}\subset \cF_{\alpha}\sqcup\cF_{\beta}.
\end{equation}

\end{rem}

\begin{pro}\label{P:5.2} Assume that $\alpha<\beta<\omega_1$ and  that $\gamma<\omega_1$
\begin{enumerate} 
\item[a)] There is an $m\in\N$ so that
 $ \cF_\alpha\cap[\{m,m+1,\ldots\}]^{<\infty}\sqcup\cF_\gamma
         \subset       \cF_\beta\cap[\{m,m+1,\ldots\}]^{<\infty}\sqcup\cF_\gamma$
\item[b)]  For any $N\infsubset\N$ it follows that
$\I(\cF_\alpha\sqcup\cF_\gamma,N)<\I(\cF_\beta\sqcup\cF_\gamma,N)$
\end{enumerate}
\end{pro}

\begin{proof}
(a) follows immediately from Proposition  \ref{P:3.1}(c). To prove  (b) define
$\delta=$\\ $\I(\cF_\alpha\sqcup\cF_\gamma,N)$. Let $M\infsubset N$.
 By Proposition \ref{P:4.2} (last  part) there is an $L\infsubset M$ so that
$\cF^L_\delta\subset\cF_\alpha\sqcup\cF_\gamma$.
Then note that
\begin{align*}
\cF^L_{\delta+1}&=\{\{\ell\}\cup D: D\in\cF_\delta^L\text{ and } \ell\in L\}\\
&\subset \{ \{\ell\}\cup A\cup G: A\in\cF_\alpha\cap[L]^{<\infty} ,
 G\in\cF_\gamma\cap[L]^{<\infty}    \ell\in L    \}\\
&\subset\{\tilde A\cup G:\tilde A\in\cF_{\alpha+1}\cap[L]^{<\infty},G\in\cF_\gamma\}
\subset \cF_{\alpha+1}\sqcup \cF_{\gamma}.
\end{align*}
Thus, by (a), $\delta=\I(\cF_\alpha\sqcup\cF_\gamma,N)<\delta+1\le
 \I(\cF_{\alpha+1}\sqcup\cF_\gamma,N)\le  \I(\cF_{\beta}\sqcup\cF_\gamma,N)$
which finishes the proof.\end{proof}

\begin{pro}\label{P:5.2b}
Assume that $(\cF_\alpha)_{\alpha<\omega_1}$ and $(\cG_\alpha)_{\alpha<\omega_1}$
 are two transfinite families and let $N\infsubset \N$  and $\alpha<\beta$ be such that 
 $\I(\cF_\alpha\sqcup\cF_\beta,N)$ and $\I(\cG_\alpha\sqcup\cG_\beta,N)$ 
 are stabilized in the sense of Proposition \ref{P:4.6}, i.e.
 $\I(\cF_\alpha\sqcup\cF_\beta,N)=\I(\cF_\alpha\sqcup\cF_\beta,\tilde N)$
 and $\I(\cG_\alpha\sqcup\cG_\beta,N)=\I(\cG_\alpha\sqcup\cG_\beta,\tilde N)$
 for all $\tilde N\subset N$.

Then it follows that  $\I(\cF_\alpha\sqcup\cF_\beta,N)=\I(\cG_\alpha\sqcup\cG_\beta,N)$.

\end{pro}

\begin{proof} Write $N$ as $N=\{n_1,n_2,\ldots\}$, with  $n_i\nearrow\infty$.
By Proposition \ref{P:3.8}  we can find an $M=\{m_1,m_2,\ldots\}\subset N$,
 $m_i\nearrow\infty$, so that 
 $\cG_\alpha^M\subset\cF_\alpha$ and   $\cG_\beta^M\subset \cF_\beta$.
 Define $L=\{ m_{n_i}: i\in\N\}$. Then it follows
that
 $\I(\cG_\alpha\sqcup\cG_\beta, N)=\I(\cG_\alpha^M\sqcup\cG_\beta^M, L)\le
 \I(\cF_\alpha\sqcup\cF_\beta, L)=\I(\cF_\alpha\sqcup\cF_\beta, N)$,
which finishes the proof by symmetry.
\end{proof}

\begin{pro}\label{P:5.2c} (Cancellation Lemma)

Let $\alpha,\beta<\omega_1$ and consider a map
$\Psi:\cF_\alpha\sqcup\cF_\beta\to[\N]^{<\infty}$ with the following property:

\noindent
There exists a hereditary $\cB\subset[\N]^{<\infty}$ and an  $\displaystyle N\infsubset\N$             
           so that $\Psi(\cF_\alpha\sqcup\cF_\beta)\!\subset\!\cB$ and $\I(\cB,N)\!\le\!\beta$.

If $\cC\subset[\N]^{<\infty}$ is hereditary and contains
 the set $\{A\setminus\Psi(A): A\in \cF_\alpha\sqcup\cF_\beta\}$ then there is 
an $M\infsubset N$ so that
 $\I(\cC,M)\ge \alpha$.
\end{pro}

\begin{proof}
Assume that  our claim is not true and that for all $M\infsubset N$ it follows
 that $\I(\cC,M)< \alpha$. First, by applying Proposition \ref{P:4.6} and passing
 to a subsequence of $N$, we can assume that
 $\I(\cB,\tilde N)=\I(\cB,N)\le \beta$ for all $\tilde N\infsubset N$. 
By applying  Proposition \ref{P:4.6}
 a second time we can also assume that 
  $\alpha_0=\I(\cC,\tilde N)=\I(\cC,N)<\alpha$, for all $\tilde N\infsubset N$, and
 define $\beta_0=\I(\cB,N)\le\beta$.

Using Theorem \ref{T:4.1} (since we are in the stabilized situation we can apply
it simultaneously to $\cC$ and $\cB$ as in the proof of Corollary \ref{C:4.10})
we obtain a transfinite family $(\cG_\gamma)$ and
a $M\infsubset N$ so that
$\cB\cap [M]^{<\infty}\subset\cG_{\beta_0}$ and
 $\cC\cap [M]^{<\infty}\subset\cG_{\alpha_0}$ and thus
$(\cB\sqcup\cC)\cap [M]^{<\infty}\subset\cG_{\beta_0}\sqcup \cG_{\alpha_0}$.

On the other hand it is clear that 
  $\cF_\alpha\sqcup\cF_\beta\subset\cC\sqcup\cB$.
 Now, using Proposition \ref{P:4.6}
 we  first find  a subset $L\infsubset M$ so that
$\I(\cG_{\beta_0}\sqcup \cG_{\alpha_0},\tilde L)=I(\cG_{\beta_0}\sqcup \cG_{\alpha_0},L)
 =\I(\cG_{\beta_0}\sqcup \cG_{\alpha_0}, M)$ for all $\tilde L\subset L$.
 Then we pass to a subset $K\infsubset L$ so that
 and $\I(\cF_\alpha\sqcup\cF_\beta,\tilde K)=\I(\cF_\alpha\sqcup\cF_\beta,K)=
 \I(\cF_\alpha\sqcup\cF_\beta, L)$, for all $\tilde K\infsubset K$.

Finally we deduce the following chain of inequalities 
\begin{equation*}
 \I(\cF_\alpha\sqcup\cF_\beta,K)\le \I(\cC\sqcup\cB,K)\le 
 \I(\cG_{\alpha_0} \sqcup\cG_{\beta_0},K)=
 \I(\cF_{\alpha_0} \sqcup\cF_{\beta_0},K)
\end{equation*}
(the last inequality follows from Proposition \ref{P:5.2b}).
But on the other it follows from Proposition \ref{P:5.2} that
 $\I(\cF_\alpha\sqcup\cF_\beta,K)>\I(\cF_{\alpha_0} \sqcup\cF_{\beta},K)
\ge \I(\cF_{\alpha_0} \sqcup\cF_{\beta_0},K)$ which is a contradiction.
\end{proof}

\section{Proof of Theorem \ref{T:1.2}}\label{S:6}

We now turn to the proof of Theorem \ref{T:1.2}. We will first restate it in an equivalent form
 using the equivalence of the Cantor Bendixson index and the concept of 
$\alpha$-largeness. Recall the definitions of   $\cB(\zb,b)$, $B(X,\beta)$, $\beta_0(\zb)$ and $\beta_0(X)$
in Definition \ref{D:1.1a} for $X$ being a Banach space, $\zb=(z_n)\subset X$
 seminormalized, $b>0$ and $\beta<\omega_1$. Using Corollary \ref{C:5.6} we get 
\begin{equation}
\label{E:6.1} b(\zb,\beta)=
\inf\Bigl\{ b\ge 1: \exists N\infsubset \N\quad \cB(\zb,b)\text{ is $\beta$-large on $N$}\Bigr\},
\end{equation}
where for $b\ge 1$ and $\cB(\zb,b)=\{A\in [\N]^{<\infty}: (z_n)_{n\in A}\sim_b\ell_1^A-\text{unit vector basis}\}$.

We secondly want to replace in the definition of $b(\zb,\beta)$ the set $\cB(\zb,b)$ by
a somewhat more convenient set.  We will need the following special case of a result 
 from \cite{AMT} (see also \cite{AG} Lemma 3.2).

\begin{lem}\label{L:6.1} 
Assume $(x_n)$ is a weakly null and  semi-normalized  sequence in a 
 Banach space $X$. Let $\delta>0$, $N\infsubset\N$ and $\vp_n>0$, for $n\in\N$.
Then there exists an $M\infsubset \N$ , $M=\{m_1,m_2,\ldots\}$ so that
 for all finite $F\subset M$ the following implication is true:
\begin{align}\label{E:6.1.1} &\text{If there is an }
 x^*\in B_{X^*}, \text{ with } x^*(x_m)\ge \delta  \text{ for all }m\in F,
 \text{then there is a }y^* \in B_{X^*},\\
 &\text{ with } y^*(x_m)\ge
\delta  \text{ for all } m\in F\text{ and } 
  |y^*(x_{m_i})|<\vp_i
 \text{ for all }i\in\N, \text{ with } m_i\not\in F. \notag
\end{align}
\end{lem}
For the sake of being self contained and reader-friendly we 
present a proof of this special case of  the above cited result  of
 \cite{AG}.
\begin{proof}
By recursion we choose for every $k\in\N$, $m_k\in\N$ and $L_k\infsubset N$
with $m_1<m_2<\ldots m_k$, 
$L_k\infsubset L_{k-1}\ldots\infsubset L_1\infsubset L_0= N$, and $m_k=\min
L_{k-1}$ so that  for all $F\subset \{1,\ldots,k\}$  all 
$L\infsubset L_k$, $L=\{\ell_1,\ell_2,\ldots\}$, and all $n\in\N$  the
following  implication holds:
\begin{align}\label{E:6.1.2}
&\left[ \exists x*\in B_{X^*} \quad x^*(x_{m_i})\ge \delta  \text{ for all }i\in F
 \text{ and }x^*(x_{\ell_i})\ge \delta \text{ for all } i=2,3,\ldots n\right]
\\
    \Rightarrow &\Biggl[ \exists y^* \in B_{X^*} \quad
y^*(x_{m_i})\ge \delta  \text{ for all } i\in F,\quad
y^*(x_{\ell_i})\ge \delta \text{ for all } i=2,3,\ldots n\notag\\
&\qquad|y^*(x_{m_i})|<\vp_i  \text{ for all }i\in\{1,\ldots k\}\setminus F
\text{ and }|y^*(x_{\ell_1})|<\vp_{k+1}\Biggr].\notag
\end{align}

Clearly the claim of the Lemma follows if we can accomplish such a choice
of $m_k$'s and $L_k$'s.

Assume for some $k\ge 1$ we have chosen $L_{k-1}$ (recall: $L_0=N$) and
 $m_1<\ldots m_{k-1}$.

We define $m_k=\min L_{k-1}$ and
\begin{equation*}
\cL=\left\{ L\subset L_{k-1}: \forall
F\subset\{1,\ldots k\}\quad \forall n\in\N \quad (L,F,n)  
\text{ satisfies (\ref{E:6.1.2}) } \right\}. 
\end{equation*}

It is easy to see that $\cL$ is closed in the pointwise topology and 
we can apply Ramsey's theorem.

In the case that there is an $L\in\cL$ so that $[L]^{\infty}\subset\cL$ we are
done. We have to show that the alternative in Ramsey's theorem leads 
to a contradiction.

 Assume that there is an $\tilde L\infsubset L_{k-1}$  so that
 $[\tilde L]^{\infty}\cap\cL=\emptyset$. Thus for
any $L=\{\ell_1,\ell_2\ldots\}\infsubset\tilde L$ there is 
 an $F=F_L\subset\{1,2\ldots k\}$ and an $n=n_L\in\N$ so that
 there exists an $x^*=x^*_L\in B_{X^*}$ with
\begin{equation}\label{E:6.1.3}
x^*(x_{m_i})>\delta \text{ for all }i\in F,
\text{ and } x^*(x_{\ell_i})>\delta  \text{ for all }i=2,3,\ldots n
\end{equation}
but for any $y^*\in B_{X^*}$ satisfying (\ref{E:6.1.3}) there must
be either an $i\in\{1,\ldots k\}\subset F$ with $|y^*(x_{m_i})|>\vp_i$ or
 $|y^*(x_{\ell_1})|\ge \vp_{k+1}$.

We first use again Ramsey's theorem to assume without
loss of generality that  the sets $F_L$ do not depend
 on $L$. Thus $F_L=F$ for all $L\infsubset \tilde L$.

Fixing for a moment such an $L=\{\ell_1,\ell_2\ldots\}\infsubset\tilde L$
 we let $j_0=\max\{j\in[0, k]: j\not\in F\}$. If $j_0=0$
 (meaning $F=\{1,2\ldots k\}$) we put $z^*_L=x^*_L$ and  observe
that we must have $|z^*_L(x_{\ell_1})|\ge \vp_{k+1}$. If $j_0\ge 1$ we apply
 the fact that our  induction hypothesis is true for $j_0-1$ and are
 able  to find a $z_L^*\in B_{X^*}$  satisfying 
 (\ref{E:6.1.3}) 
 and secondly 
\begin{equation}\label{E:6.1.4}
|z_L^*(x_{m_i})|<\vp_i \text{  if }i\in\{1,2,\ldots k\}\setminus F
\end{equation} 
and,  thus, we also   must have  $|z_L^*(x_{\ell_1})|>\vp_{k+1}$
(apply the induction hypothesis to the set\\
 $L=\{m_{j_0},\ldots m_k,\ell_2,\ell_3, \ldots\}\subset L_{j_0-1}$,
  $F\cap [1,j_0-1]$ and $n+k-j_0+1$).

We write $\tilde L=\{\tilde\ell_1,\tilde\ell_2,\ldots\}$ and claim
that for  any $m\in\N$ there is a $x_m^*$ so that
$|x^*_m(x_{\tilde\ell_i})|\ge \vp_{k+1}$,
 for $i=1,2,\ldots m$. This would be a contradiction
to the assumption that $(x_i)$ is weakly null.
                                                                                                          
For each $j=1,2,\ldots m$ we find an $n_j$ so that 
 the triple $(\tilde L_j,F,n_j)$ does not satisfy (\ref{E:6.1.2})
with $\tilde L_j=\{\tilde\ell_j,
\tilde\ell_{m+1},\tilde\ell_{m+2}\ldots \}$. Choose $j_0$ so that
$n_{j_0}$ is the maximum of $(n_j)_{j=1}^m$. Then choose
 $x^*_m=z^*_{L_{j_0}}$ (where $z_L$ is defined as above).
                                                                                                          
We observe that  $x^*_m$ satisfies   
 (\ref{E:6.1.3}) with respect to all of the $L_j$'s. Secondly it
satisfies (\ref{E:6.1.4}) and thus  it must  follow that
 $|x^*_m(x_{\tilde\ell_j})|>\vp_{k+1}$ for all $j=,2,\ldots m$.
This finishes the proof of the claim and, thus the proof of the
Lemma.\end{proof}

We will slightly reformulate Lemma \ref{L:6.1}. 
\begin{cor}\label{C:6.1a}
Assume $(x_n)$ is a weakly null and  semi-normalized  sequence in a 
 Banach space $X$. Let $\delta>0$, $N\infsubset\N$ and $\vp>0$.

Then there exists an $M\infsubset \N$ , $M=\{m_1,m_2,\ldots\}$ so that
 for all finite $F\subset M$  the following implication is true:
\begin{align}\label{E:6.1a.1} &\text{If there is an }
 x^*\!\in\! B_{X^*}, \text{ with } x^*(x_m)\!\ge \!\delta  \text{ for all }m\!\in\! F,
  \text{then there is also a }  z^*\!\in\!B_{X^*},\\
& \text{ with } z^*(x_m)\ge
\delta-\vp  \text{ for all } m\in F\text{ and } 
   z^*(x_{m_i})=0,
 \text{ for all }i\in\N,\text{ with } m_i\not\in F. \notag
\end{align}
\end{cor}
\begin{proof} After passing to a subsequence of $x_n$ we can assume that
there is for each $n\in\N$ an $x^*_n\in (2/\|x_n\|)B_{X^*}$, with $x^*_n(x_m)=\delta_{(n,m)}$, whenever
 $n,m\in\N$. 
 Choose for $n\in\N$ $\vp_n=2^{-n}\vp/(2+\sup_{n\in\N}\|x_n\|)$ and apply Lemma \ref{L:6.1} in order to obtain 
$M=(m_i)\infsubset\N$. 
 
If $F\subset M$ is finite and $x^*\in B_{X^*}$ so that $x^*(x_i)\ge \delta$, whenever $i\in F$, then we
 let $y^*\in B_{X^*}$ be as prescribed in (\ref{E:6.1.1}) and let 
 $z^*=\tilde z^*/\|\tilde z^*\|$ where
$\tilde z^*=y^*-\sum_{m\in M\setminus F} y^*(x_{m})x^*_m$.
\end{proof}
                      
We introduce notations similar to $\cB(\zb,b)$, $b(\zb,\beta)$, $B(\beta,X)$, $\beta_0(\zb)$
 and $\beta_0(X)$.
\begin{defn}\label{D:6.2}
Let $\xb=(x_n)$ be a seminormalized sequence in a Banach space $X$.
For $a>0$ we put 
\begin{equation} \label{E:6.2.1}
\cA(\xb,a)=
\bigl\{A\in [\N]^{<\infty}:\exists x^*\in B_{X^*}\forall i\in A \quad x^*(x_i)\ge a\bigr\}. 
\end{equation}
For $\alpha<\omega_1$  we let
\begin{align}\label{E:6.2.2}
&a(\xb,\alpha)=\sup\bigl\{ a\ge 0: \exists N\infsubset\N \quad\cA(\xb,a)
\text{ is $\alpha$-large on $N$}\bigr\} \\
         &\qquad\qquad=\sup\bigl\{ a\ge 0: \exists M\infsubset \N \quad\cF_\alpha^M\subset\cA(\xb,a)\bigr\}, 
      \notag\\
\label{E:6.2.3}&\alpha_0(\xb)=\sup\{\alpha<\omega_1: a(\xb,\alpha)>0\}, \\
\label{E:6.2.4} &A(\alpha,X)=\sup \{ a(\zb,\alpha): \zb\subset B_X \text{ seminormalized and weakly null}\},\text{ and}\\
\label{E:6.2.5}
 &\alpha_0(X)=\sup\{\alpha\!<\!\omega_1: A(\alpha,X)>0\}=\sup\{\alpha_0(\zb):
   \zb\!\subset\! B_X \text{ seminorm., weakly null}\}.
\end{align}
\end{defn}

\begin{lem}\label{L:6.3} Let $\xb$ be a seminormalized sequence in a Banach space $X$ with
 $\alpha_0(\xb)<\omega_1$.

Then there is subsequence $\yb=(y_n)$ of $\xb$  with the following properties.
\begin{enumerate}
\item[a)] For all $\alpha<\omega_1$ and all $a'<a(\yb,\alpha)$ the set $\cA(\yb,a')$ is $\alpha$-large on
 $\N$.
\item[b)] For all $\alpha<\omega_1$ and all subsequences $\zb$ of $\yb$ it follows that
 $a(\yb,\alpha)=a(\zb,\alpha)$.
\item[c)] The map $[0,\omega_1)\ni\alpha\mapsto a(\yb,\alpha)$ is decreasing and continuous.
\end{enumerate}
Moreover if $\beta_0<\alpha_0(\xb)$ (and, thus $a(\xb,\beta_0)>0$) and if $0<\eta<a(\xb,\beta_0)$ 
 the subsequence $\yb$ can be chosen so that $a(\yb,\beta_0)>a(\xb,\beta_0)-\eta$.

\end{lem}
\begin{proof} We first note that if $(\ub)=(u_n)$ is almost a subsequence of $\vb=(v_n)$
 (i.e. for some $n_0\in\N$ it follows that $(u_{n_0+i})_{i\in\N}$ is a subsequence of $\vb$)
 then $a(\ub,\alpha)\le a(\vb,\alpha)$. In particular this means that
$\alpha_0(\ub)\le \alpha_0(\vb)$ and therefore it is enough to 
 find a subsequence of $\xb$ which satisfies (a), (b) and (c) for all
 $\alpha<\alpha_0(\xb)$.

\noindent Claim.  Let $\yb=(y_n)$ be a subsequence of $\xb$ and let $\alpha<\alpha_0(\xb)$. 
Then there is a subsequence $\zb$ of $\yb$ so that for all $a'<a(\zb,\alpha)$ it follows that
 $\cA(\zb,a')$ is $\alpha$-large on $\N$.

Because of the observation at the beginning of the proof the sequence $\zb$ in the claim has
 the property that  $a(\ub,\alpha)=a(\zb,\alpha)$ for any sequence $\ub$ which is
almost a subsequence of $\zb$.

In order to show the claim we let $\vp_i\searrow 0$, put $a_0=a(\yb,\alpha)$ and choose
 an $N_1=(n^{(1)}_i)\infsubset \N$  so that $\cF^{N_1}_\alpha\subset \cA(\yb,a_0-\vp_1)$.
 Letting now $\zb^{(1)}=(y_{n^{(1)}_i})$ we deduce that
 $\cF_\alpha\subset \cA(\zb^{(1)},a_0-\vp_1)$ and, thus,
 $\cF_\alpha\subset\cA(\zb,a_0-\vp_1)$ for any subsequence $\zb$ of $\zb^{(1)}$
 (recall that $\cF_\alpha$ is spreading), which
finally implies that $\cA(\zb,a_0-\vp_1)$ is $\alpha$-large for any sequence $\zb$ 
which is almost a subsequence of $\zb^{(1)}$. 
 
Now we let $a_1=a(\zb^{(1)},\alpha)$ and continue this way, eventually finding
 $\N\infsupset N_1=(n^{(1)}_i)\infsupset N_2$ $= (n^{(2)}_i)\infsupset\ldots$,
so that if we put $\zb^{(k)}=(y_{n^{(k)}_i})$ 
 and $a_k= a(\zb^{(k)},\alpha)$ it follows that
$\cF_\alpha\subset \cA(\zb^{(k)}, a_{k-1}-\vp_k)$. We deduce that $a_{k-1}-\vp_k\le a_k\le a_{k-1}$. 
Letting $\zb$ be a diagonal sequence of the $\zb^{(k)}$'s it follows that
 $a(\zb,\alpha)=a=\inf a_k$ and, since every subsequence of $\zb$ is almost a subsequence of
  each $\zb^{(k)}$ it follows for each $k\in \N$ that  
$\cA(\zb,a_{k-1}-\vp_k)$ is $\alpha$-large,  which implies the claim.

Note also, that if we had assumed that none of  the $\vp_i$'s would
 exceed a value $\eta$ then it follows that $a(\zb,\alpha)>a(\yb,\alpha)-\eta$
 (this proves  the part of our claim starting with ``moreover'' if we let $\alpha=\beta_0$).
 Writing now the interval $[0,\alpha_0(\xb))$ as a sequence $(\alpha_n)$ 
  and applying
successively the above claim to each $\alpha_n$, we obtain by diagonalization
 a subsequence $\yb$ of $\xb$ so that  (a) and (b) of our statement are satisfied.
 It is also clear that $a(\yb,\alpha)$ is decreasing in $\alpha$. Let $\alpha$ be a limit ordinal
and $a'<a=\lim_{\beta\to \alpha} a(\yb,\beta)=\inf_{\beta<\alpha} a(\yb,\beta)$,
then it follows that for every $\beta<\alpha$ that 
$\cA(\yb,a')$ is $\beta$-large on $\N$. By Proposition \ref{P:4.2} this implies that
 $\cA(\yb,a')$ is $\alpha$-large on $\N$ for all $a'<a$, which implies the claimed continuity.
\end{proof}

\begin{pro}\label{P:6.2a} Let $\xb=(x_n)$ be a weakly null, and normalized sequence in a Banach space $X$
 and $c>\eta>0$.
Then there is a subsequence $\yb$ of $\xb$ so that 
\begin{equation}\label{E:6.2a.1}
\cA\bigl(\yb,2c+\eta\bigr)\subset\cB \bigl(\yb,\frac1c\bigr)\subset\cA\bigl(\yb,c-\eta\bigr)
\end{equation}
and, thus it follows that 
\begin{align}\label{E:6.2a.2}
&\frac12 a(\xb,\alpha)\le \frac1{b(\xb,\alpha)}\le a(\xb,\alpha)\text{ and }
 \frac12 A(X,\alpha)\le \frac1{B(X,\alpha)}\le A(X,\alpha),
\\
\label{E:6.2a.3}&\alpha_0(\xb)=\beta_0(\xb),\text{ and }\alpha_0(X)=\beta_0(X)
\end{align}
\end{pro}

Proposition \ref{P:6.2a}  will follow from Corollary \ref{C:6.1a} and the following
simple observation.

\begin{lem}\label{L:6.3b}
Assume $E=(\R^n,\|\cdot\|)$, $n\in\N$ is an $n$-dimensional normed space 
for which the unit vector basis $(e_i)_{i=1}^n$ of $\R^n$ is a normalized basis.
 Define:
\begin{align*}
&c_1=\max\Biggl\{ c\ge 0: \forall A\subset\{1,2\ldots n\}:\exists x^*\in B_{E^*}
 \quad \begin{matrix} &x^*(e_i)\ge c \text{ if } i\in A\\
                      & x^*(e_i)=0 \text{ if } i\not\in A\end{matrix}\Biggr\}\\
&c_2=\max\bigl\{ c\ge 0:cB_{\ell_\infty^n}\subset B_{E^*}\bigr\}\text{ and }
   c_3=\min\Bigl\{ \|\sum_{i=1}^na_ie_i\|: \sum_{i=1}^n|a_i|=1\Bigr\}
\end{align*}
Then it follows that $c_1\ge c_2=c_3\ge \frac12 c_1$.
\end{lem}
\begin{proof} It is clear that $c_1\ge c_2$. 
 To show that $c_2\ge c_3$   we first observe  that by the maximality
 of $c_2$  we can  find an $x^*\in S_{X^*}$ of the form
 $x^*=\sum_{i=1}^n x^*_ie_i^*\in S_{E^*}$ 
($e_i^*$ being the $i$-th coordinate functional) so that $|x^*_i|=c_2$ for $i=1,\ldots n$.
 Then choose $x=(x_i)\in S_E$ so that $x^*(x)=\sum x_i x^*_i=1$. Thus
 $1=\sum x_i x^*_i\le c_2\sum |x_i|\le c_2/c_3$,
which implies the claimed inequality.

In order to show $c_3\ge c_2$  and $c_3\ge  c_1/2$ let $(a_i)_{i=1}^n\in \R^n$ with
 $\sum_{i=1}^n |a_i|=1$.
First we can choose an $x^*=\sum_{i=1}^n x^*_i e^*_i\in B_{E^*}$ with 
 $ | x^*_i|=c_2$ and $\sign(x^*_i)=\sign(a_i)$, for $1\le i\le n$. This proves
that $\|\sum_{i=1}^n a_ie_i\|\ge x^*(\sum_{i=1}^n a_ie_i)\ge c_2$, which implies $c_3\ge c_2$. Secondly, we can
 assume that  $\sum_{i=1}^n a_i^+\ge 1/2$ and a simliar argument
 implies that $\|\sum_{i=1}^n a_ie_i\|\ge c_1/2$, and thus $c_3\ge  c_1/2$.
\end{proof}

\begin{proof}[Proof of Proposition \ref{P:6.2a}] Assume $c>\eta>0$ to be given.
Using  Corollary \ref{C:6.1a}  we can find a subsequence $\yb=(y_n)$ of $\xb$ so that
\begin{equation*}
\cA(\yb,2c+\eta)\subset\tilde\cA(\yb,2c),\text { and (trivially) }
\tilde\cA(\yb,c)\subset\cA(\yb,c),\text{ where we put for $r>0$}
\end{equation*}
 $\tilde\cA(\yb,r)\!=\!\{A\!\in\![\N]^{<\infty}\!: \forall B\subset A\exists x^*\!\in\! B_{X^*}\quad
 x^*(y_i)\!\ge \!r\text{ if }i\!\in \!B\text{ and }x^*(y_i)\!=\!0 \text{ if }i\!\not\in \!B \}$.
 
 Secondly we deduce from Lemma \ref{L:6.3b} that
\begin{equation*}
\tilde\cA(\yb,2c)\subset\{A\in[N]^{\infty}:\forall(a_i)_{i\in A}\subset \R\quad
       \|\sum_{i\in A} a_iy_i\|\ge c\sum_{i\in A} |a_i|\}\subset\tilde\cA( \yb,c),
\end{equation*}
proving the claim (note that $1/c_3$, where $c_3$ as defined in Lemma \ref{L:6.3b}
  is the smallest $c$ so that $E$ is $c$-equivalent to $\ell^n_1$).
\end{proof}

We now can restate Theorem \ref{T:1.2} as follows.

\begin{thm}\label{T:6.4} Let $X$ be a  Banach space with a basis not containing $c_0$. Assume that 
there is an  ordinal $\beta_0\in [0,\alpha_0(X)]$ so that
the following two conditions hold.
\begin{enumerate}
\item[a)] There is a seminormalized weakly null sequence $\yb\subset B_X$ with
  $a(\yb,\beta_0)=0$.
\item[b)] $\inf_{\gamma<\beta_0} B(X,\gamma)>0$.
\end{enumerate}
Then there is a seminormalized block basis $(x_n)$ in $X$
 and a subsequence $(\tilde y_n)$ of $(y_n)$ so that
the map $x_n\mapsto \tilde y_n$ extends to a strictly singular and linear bounded operator
 \begin{equation*}T: [x_n:n\in\N]\to [\tilde y_n:n\in\N].\end{equation*}
\end{thm}

In order to prove Theorem \ref{T:6.4} we will need several Lemmas.

\begin{lem}\label{L:6.5}
Let $(\cF_\alpha)$ be a transfinite family and $\xb=(x_n)$ be a weakly null 
and semi-normalized sequence in a Banach space $X$
 satisfying the conclusions of Lemma  \ref{L:6.3}.
Let 
$1\le\alpha<\omega_1$ and  
   assume that $a(\xb,\alpha)>0$.
For any $\eta>0$ there is a subsequence $\zb$ of $\xb$ so that:

 For  $A\in \cF_\alpha$, or $A\in\cF_\alpha\sqcup\cF_\alpha$,  there
 is a $z^*_A\in B_{X^*}$ so that $z^*_A(x_m)=0$ if $m\in \N\setminus A$ and
 $z^*_A(x_m)\ge a(\xb,\alpha)-\eta$, or
 $z^*_A(x_m)\ge (a(\xb,\alpha)-\eta)/2$, if $m\in A$, respectively.  
\end{lem}
\begin{proof} 
Let $1\le\alpha<\omega_1$, with $a(\xb,\alpha)>0$, and $\eta>0$. 
 Put $a=a(\xb,\alpha)$. From the definition of $a(\cdot,\cdot)$  it follows that
there is an $M_1\subset\N$ so that
$\cF_\alpha^{M_1}\subset \cA(\xb,a-\eta/4)$
 and applying Corollary  \ref{C:6.1a} 
 we find an $M_2\infsubset M_1$ so that
 (\ref{E:6.1a.1}) holds for $\delta=a-\eta/4$ and $\vp=\eta/4$. In particular this
 means that  our claim holds for all $A\in\cF_\alpha^{M_2}$.
 Secondly we deduce that
 $\cF_\alpha^{M_2}\sqcup\cF_\alpha^{M_2}\subset \cA(\xb,(a-\eta/2)/2)$. Since 
each $A\in\cF_\alpha\sqcup\cF_\alpha$ is the disjoint union of two elements of $\cF_\alpha$,
 it also follows that for any $A\in\cF_\alpha^{M_2}\sqcup\cF_\alpha^{M_2}$
 we find a $z^*_A\in B_{X^*}$ so that $z^*_A(x_m)=0$ if $m\in M_2\setminus A$
  and $z^*_A(x_m)\ge (a(\xb,\alpha)-\eta)/2$, if $m\in A$.

Choosing finally $\zb$ to be the subsequence defined by $M_3$, will finish the proof. 
 \end{proof}
Using  Lemma \ref{L:6.5} successively for different $\alpha$'s  and the appropriate choices of $\eta$
 we conclude from  a simple diagonalization argument the following Corollary. 
\begin{cor}\label{C:6.5a} Let $\xb=(x_n)$ be a weakly null and semi-normalized sequence in a Banach space $X$
 satisfying the conclusions of Lemma  \ref{L:6.3}, let $(\ell_k)\subset\N$ be strictly increasing,
and let $(\alpha_k)\subset [0,\alpha_0(\xb))$. Then there is a subsequence $\zb$ so that
 for any $k\in\N$ and any $\ell\le\ell_k$ it follows that: 

For  $A\in \cF_{\alpha_\ell}\cap [\{k,k+1,\ldots\}]^{<\infty}$, 
or $A\in[\cF_{\alpha_\ell}\sqcup\cF_{\alpha_\ell}]\cap [\{k,k+1,\ldots\}]^{<\infty}$,  there
 is a $z^*_A\in B_{X^*}$ so that $z^*_A(x_m)=0$ if $m\in \N\setminus A$ and
 $z^*_A(x_m)\ge a(\zb,{\alpha_\ell})/2$, or
 $z^*_A(x_m)\ge a(\zb,{\alpha_\ell})/4$, if $m\in A$, respectively.  

\end{cor}
\begin{lem}\label{L:6.6}
Let $\xb=(x_n)$ be a weakly null and semi-normalized sequence in a Banach space $X$
 satisfying the conclusions of Lemma  \ref{L:6.3}. Let $(\delta_k)$ be decreasing sequence 
in $(0,1)$, with $\delta_k\le 1/(k+1)$, for any $k\in\N$.

Then there is a subsequence $\zb$ of $\xb$
 and a sequence
of ordinals $\alpha_k$ which increases to $\alpha_0(\xb)$
 so that 
for each $k\in\N$
\begin{enumerate}
\item[a)] $\frac12 \delta_k \le a(\zb,\alpha_k)  \le \delta_k$
\item[b)] $\alpha_k$ is of the form $\alpha_k=\tilde \alpha_k + k$ 
 (here we identify positive integers with finite ordinals of the same cardinality).
\end{enumerate}
\end{lem} 
\begin{proof} 
We define for $k\in\N$  
$\beta_k=\min\{\beta : a(\xb,\beta)\le \delta_k\}$.
Since $\delta_k\le 1/(k+1)$ it is easy to see that $\beta_k\ge (k+1)$.
If $\beta_k$, is a successor,  say $\beta_k=\gamma_k+1$, then 
$a(\xb,\gamma_k)\ge \delta_k$
If $\beta_k$ is a limit ordinal we deduce from the continuity  that $a(\xb,\beta_k)=\delta_k$,
 and we let $\gamma_k=\beta_k$.

It follows from (\ref{E:5.7.1}) that
 $a(\xb,2\gamma_k)\ge \frac12a(\xb,\gamma_k)\ge \frac12\delta_k$.
 Since $\gamma_k\ge k$ we will be able to find an $\alpha_k$ between $\beta_k$ and $2\gamma_k$
 which is of the form as required in (b).
\end{proof}

\begin{lem}\label{L:6.7} Assume that $\xb=(x_n)$ is a weakly null and seminormalized sequence 
 in a Banach space $X$ satisfying the conclusion of Lemma \ref{L:6.3}.
Let $\alpha_k$ be  ordinals increasing to $\alpha_0(\xb)$
  so that for each $k\in\N$ $\alpha_k$ can be written as $\alpha_k=\tilde\alpha_k+k$
for some $\tilde\alpha_k<\omega_1$.
 For
 $k\in\N$  let $\vp_k=a(\zb,\alpha_k)$.

\begin{enumerate}
\item[a)]
 There exists a transfinite family $(\cG_\alpha)$ and a subsequence
 $\zb=(z_n)$ of $\xb$ so that for all $(a_i)\in\coo$ it follows that
\begin{equation*}
\Bigl\|\sum a_i z_i\Bigr\|\le8 \sum_{k=1}^\infty \vp_{k-1}\max_{A\in\cG_{\alpha_k}}\sum_{i\in A} |a_i|.
 \end{equation*}
\item[b)] Let $(\cF_\alpha)$  be a transfinite family. Then there is a subsequence
 $\zb=(z_n)$ of $\xb$ so that for all $(a_i)\in\coo$ it follows that
\begin{equation*}
 \Bigl\|\sum a_i z_i\Bigr\|\ge \frac18\max_{k\in\N,A\in\cF_{\alpha_k}}\vp_k\sum_{i\in A} |a_i|.  
 \end{equation*}
\end{enumerate}
\end{lem}
\begin{proof}  From our assumption on $\xb$ it follows that for all $Q\infsubset \N$ and
 all $k\in\N$ it follows that $\I(\cA(\xb,2\vp_k),Q)\le \alpha_k$.
Therefore we can apply Corollary \ref{C:4.10} to obtain a transfinite family 
$(\cG_\alpha)_{\alpha<\omega_1}$ and a $K=\{k_1,k_2,\ldots\}\infsubset \N$, $k_i\nearrow \infty$,
 if $i\nearrow\infty$, so that for all $\ell$
\begin{equation}\label{E:6.7.3} 
\cA(\xb,2\vp_\ell)\cap[\{k_\ell,k_{\ell+1},\ldots\}]^{<\infty}\subset \cG^K_{\alpha_\ell}.
\end {equation}
 Putting $\zb=(x_{k_i})$ this implies that for all $\ell$
\begin{equation}\label{E:6.7.4} 
\cA(\zb,2\vp_\ell)\cap[\{\ell,{\ell+1},\ldots\}]^{<\infty}\subset \cG_{\alpha_\ell}.
\end {equation}
Let $(a_i)\!\in\!\coo$. We find a $z^*\!\in\! B_{X^*}$ which norms
 $\sum a_i z_i$ and ($\vp_0\!=\!\sup\|z_i\|$) and deduce that
\begin{align}\label{E:6.7.5} 
z^*\Bigl(\sum_{i=1}^\infty a_i  z_i\Bigr)
&=\sum_{k=1} \sum_{2\vp_k<z^*(z_i)\le 2\vp_{k-1}}z^*(z_i) a_i+ 
             \sum_{2\vp_k<-z^*(z_i)\le 2\vp_{k-1}} z^*(z_i) a_i\\
&\le\sum_{k=1}^\infty 2\vp_{k-1}\Bigl[\sum_{2\vp_k<z^*(z_i)\le 2\vp_{k-1}} |a_i|+
           \sum_{2\vp_k<-z^*(z_i)\le 2\vp_{k-1}} |a_i|\Bigr]. \notag
\end{align}
Now note that the set $\{i\in\N: 2\vp_k<z^*(z_i)\}$ is the union
 of two sets in $\cG_{\alpha_k}$ namely
 $\{i\ge k: 2\vp_k<z^*(z_i)\}$ and  $\{i\le k: 2\vp_k<z^*(z_i)\}$ (where the second
set lies in $\cG_{\alpha_k}$ because $\alpha_k$ is the $k$-th successor of some $\tilde \alpha_k$).
Similarly we  can proceed with $\{i\in\N: 2\vp_k<-z^*(z_i)\}$  and we therefore derive 
 from (\ref{E:6.7.5})  that
\begin{align}\label{E:6.7.6} 
z^*\Bigl(\sum_{i=1}^\infty a_i  z_i\Bigr)
\le \sum_{k=1}^\infty 2\vp_{k-1} 4\max_{A\in\cG_{\alpha_k}}\sum_{i\in A} |a_i|,
\end{align}
which finishes the proof of part (a).

In order to prove part (b) we first choose $\ell_0=0$ and $(\ell_k)_{k\in\N}\subset \N$
 fast enough increasing  so that  for all
 $(a_i)\in\coo$  it follows that
\begin{equation}\label{E:6.7.7}
\Bigl\|\sum a_i  z_i\Bigr\|\ge \sup_{A\subset \N, \#A=k}  \frac{\vp_{\ell_{k}}}4\sum_{i\in A} |a_i|.
\end{equation}
After passing to a subsequence we can assume that the conclusion of Corollary \ref{C:6.5a} is satisfied.

Let $(a_i)\in\coo$, $k\in\N$, $\ell\le\ell_k$ and $A\in\cF_{\alpha_\ell}\cap[\{k,k+1,\ldots \}]^{<\infty}$.
 Letting $A^+=\{i\in A: a_i\ge 0\}$ and $A^-=A\setminus A^+$ we can choose
 $z^*_{A^+}$ and $z^*_{A^-}$ as prescribed in Corollary \ref{C:6.5a} and deduce that
\begin{equation}\label{E:6.7.8}
\Bigl\|\sum_{i=k}^\infty a_i  z_i\Bigr\|\ge
\frac12 z^*_{A^+}\Bigl(\sum_{i=k}^\infty a_i  z_i\Bigr)-
\frac12 z^*_{A^-}\Bigl(\sum_{i=k}^\infty a_i  z_i\Bigr)\ge \frac14\vp_\ell\sum_{i\in A} |a_i|.
\end{equation}
Therefore we  obtain
\begin{align*}
\Bigl\|\sum_{i=1}^\infty a_i  z_i\Bigr\|&\ge
 \frac12\max_{k\in\N} \Bigl[\Bigl\|\sum_{i=k}^\infty a_i  z_i\Bigr\|+
                   \Bigl\|\sum_{i=1}^{k-1} a_i  z_i\Bigr\|\Bigr]\\
&\ge\frac18\max_{k\in\N}
\Bigl[\max_{\substack{\ell\in(\ell_{k-1},\ell_k],\\ A\in\cF_{\alpha_\ell}, A\ge k}}\vp_\ell\sum_{i\in A} |a_i|+
        \vp_{\ell_{k-1}}\sum_{i=1}^{k-1} |a_i|\Bigr]\\
&\ge\frac18\max_{k\in\N}\max_{\substack{\ell\in(\ell_{k-1},\ell_k]\\ A\in\cF_{\alpha_\ell}}}\vp_\ell\sum_{i\in A} |a_i|
=\frac18\max_{\ell\in\N,A\in\cF_{\alpha_\ell}}\vp_\ell\sum_{i\in A} |a_i|,
\text{ proving part (b).}
\end{align*}\end{proof}

\begin{lem}\label{L:6.8}
Assume that $\ub=(u_n)$ and $\vb=(v_n)$  are two
 weakly null and  seminormalized sequences in a Banach space $X$
satisfying the conclusions of Lemma \ref{L:6.3} such that
 \begin{equation}\label{E:6.8.1} a(\ub,\alpha)<\frac15 a(\vb,\alpha).
  \end{equation}
 Let $\xb=\ub+\vb=(u_n+v_n)$.   Then there is a subsequence $\zb$ of $\xb$  so  that
\begin{equation}\label{E:6.8.2}
\cF_\alpha \subset \cA\bigl(\zb,a(\vb,\alpha)/4\bigr), 
 \text{ in particular }
  a(\yb,\alpha)\ge a(\vb,\alpha)/4 
\text{ for all  subsequences $\yb$ of $\zb$}.
\end{equation}
\end{lem}
\begin{proof}  By passing to a subsequence we can assume that the conclusions 
 of Lemma \ref{L:6.5} are satisfied for  $\vb$, $\alpha$, and $\eta=\frac1{10}a(\vb,\alpha)$.
 In particular we find for
 each $A\in\cF_{\alpha}\sqcup\cF_{\alpha}$, 
a $v^*_A\in B_{X^*}$ which has the property
that $v^*_A(x_m)\ge (a(\vb,\alpha)-\eta)/2=\frac9{20}a(\vb,\alpha)$, whenever $m\in A$,
and  $v^*_A(x_m)=0$, if $m\not\in A$.

For $A\in\cF_\alpha\sqcup\cF_\alpha$ we define $\psi(A)=\{i\in A: v_A^*(u_i)<-\frac15 a(\vb,\alpha)\}$,
and note that \\
$\Psi(\cF_\alpha\sqcup\cF_\alpha)\subset\cA(\ub,\frac15 a(\vb,\alpha))$ 
and $\I(\cA(\ub,\frac15 a(\vb,\alpha)),\N)\le \alpha$ (using (\ref{E:6.8.1})).
Secondly we note that $\{A\setminus \Psi(A):A\in \cF_\alpha\sqcup\cF_\alpha\}$
 contains $\cA(\vb+\ub,\frac14a(\vb,\alpha))$. This means that we can apply
 Proposition \ref{P:5.2c} and derive that for  some $N\infsubset \N$
 we have  $\I(\cA(\vb+\ub,\frac14a(\vb,\alpha)),N)\ge \alpha$ and thus
that there is an $M\infsubset N$ so that
 $\cF_\alpha^M\subset  \cA(\vb+\ub,\frac14a(\vb,\alpha))$ which implies the claim
for the subsequence of $\vb+\ub$ defined by $M$.
\end{proof} 

\begin{lem}\label{L:6.9} 
 For $\alpha<\omega_1$,  $(\alpha_n)\subset[0,\omega_1)$  strictly increasing,
  $(\Theta_n)\subset [0,1]$ being non decreasing,
 and $(\ell_n)\subset \N$, with $\lim_{n\to\infty}\ell_n=\infty$ it follows that the Schreier spaces 
 $S(\cF_\alpha)$, $S((\Theta_n),(\cF_{\alpha_n}))$ and  $S((\Theta_n),(\cF_n),(\ell_n))$
 (recall Definition \ref{D:1.4})
  are hereditarily $c_0$, assuming that in the case of the space $S((\Theta_n),(\cF_{\alpha_n}))$
  the sequence
 $(\Theta_n)$ decreases to 0.
\end{lem}
\begin{proof}
We will proof by transfinite induction on $\alpha<\omega_1$ that 
the spaces $S(\cF_\alpha)$, $S((\Theta_n),(\cF_{\alpha_n}))$ and  $S((\Theta_n),(\cF_n),(\ell_n))$,
where $\alpha_n\nearrow\alpha$, is hereditarily $c_0$ (of course
 the claim for the second and third  space is vacuous if $\alpha$ is not  a limit ordinal).

Assume the claim is true for all $\beta<\alpha$. If $\alpha$ is a successor, say $\alpha=\beta+1$
  it is clear that  the norm on $S(\cF_\beta)$ and $S(\cF_{\beta+1})$ are equivalent, and, thus,
the claim follows from the induction hypothesis.

If $\alpha$ is a limit ordinal then $S(\cF_\alpha)$ is a special case 
 of the space $((\Theta_n),(\cF_{\alpha_n}),(\ell_n))$.
 Indeed,  we let $\Theta_n=1$, for $n\in\N$, and if $A_n(\alpha)$ is the approximating sequence  
 we write $\bigcup A_n(\alpha)$ as a strictly increasing sequence $\{\alpha_n:n\in\N\}$
 and put $\ell_n=\min\{\ell:\alpha_n\in A_\ell(\alpha)\}$.

 Therefore  we only need to consider the norms
\begin{equation*} 
\|x\|=\begin{cases} \sup_{n\in N}\sup_{F\in\cF_{\alpha_n},\ell_n\le F}\Theta_n\sum_{i\in F} |x_i| &\text {or }\\
 \sup_{n\in N}\sup_{F\in\cF_{\alpha_n}}\Theta_n\sum_{i\in F} |x_i| &\quad
 \end{cases}.
\end{equation*}
Let $X$ be the completion of $\coo$ under $\|\cdot\|$,
 and let, without of generality, $Y$ be the closed subspace spanned
 by a normalized block $(y_n)$. 
Either there is a further block $(z_n)$
 and an $n_0\in\N$  so that the norm $\|\cdot\|$ is on the
  closed subspace spanned by  $(z_n)$ equivalent to the norm $\|\cdot\|_{\alpha_{n_0}}$ on
 $S(\cF_{\alpha_{n_0}})$. Then
 our claim follows from the induction hypothesis.
Or we can choose a normalized  block $(z_n)$ in $Y$
  and increasing sequences
 $(k_n),(m_n)$ in $\N$, with $k_n<m_n<k_{n+1}$ so that for each $n\in\N$ it follows that
 (we denote for $z=\sum a_i e_i\in\coo$ and $k\in \N$
 $z\vert_{[\ell_k,\infty)}=\sum_{i\ge k} a_i e_i$)
\begin{align}\label{E:6.9.4}
 &\qquad\qquad\|z_{n}\|_{\alpha_k}<2^{-n}, \text{ whenever $k\le k_n$ }\\ 
  \label{E:6.9.6}  &\max\supp z_n\ge \ell_n \text{ (if  we assume admissibility) and } 
1=\|z_n\|=\Theta_{m_n}\|z_n\|_{\alpha_{m_n}}\\
  &\label{E:6.9.5}\qquad\qquad\left.   \begin{matrix} &\Theta_k\|z_{n}\|_{\alpha_k}   \\
                               &\text{ respectively }  \\
                         &\Theta_k\|z_{n}\big\vert_{[k,\infty)}\|_{\alpha_k}   
                          \end{matrix}\right\} <2^{-n}, \text{ whenever $k\ge k_{n+1}$ } 
\end{align}
(for (\ref{E:6.9.5}) we are using either the $(\ell_n)$-admissibility condition
 or the fact that $\Theta_n\searrow 0$, if $n\nearrow \infty$).

Finally it is easy to see that 
 (\ref{E:6.9.4}), (\ref{E:6.9.5})and (\ref{E:6.9.6}) imply that $(z_{m_n})$ is equivalent to the $c_0$-unit basis.\end{proof}

We are now ready for the proof of Theorem \ref{T:6.4}.

\begin{proof}[Proof of Theorem \ref{T:6.4}]

We choose  a seminormalized and weakly null sequence $\yb\subset B_X$ 
 with the property that $a(\yb,\beta_0)=0$ and
  (by passing to a subsequence)  assume that $\yb$ satisfies the conclusion of
 Lemma  \ref{L:6.3}. We let
\begin{equation}\label{E:6.4.1}
c=\min\Bigl\{ \frac15,\frac12\inf_{\gamma<\beta_0} B(X,\gamma)\Bigr\}
\end{equation} 
Using  first Lemma \ref{L:6.6} we find an increasing  sequence of ordinals 
 $(\alpha_k)\subset[0,\alpha_0(\yb))$  so that, after passing  to a subsequence, we 
  can assume that $\frac12 (c/2)^{3k} \le a(\yb,\alpha_k)  \le (c/2)^{3k}$ and
that $\alpha_k$ is the $k$-th successor of another ordinal, for  all $k\in\N$. 
 Using secondly Lemma \ref{L:6.7}(a)
 we can assume,  after passing again to a subsequence, that
 there is a transfinite family $(\cG_\alpha)_{\alpha<\omega_1}$ so that
 for any $(a_i)\in\coo$ it follows that
\begin{equation}\label{E:6.4.2}
\Bigl\|\sum a_i y_i\Bigr\|\le
 8 \sum_{k=1}^\infty (c/2)^{3k-3}\max_{A\in\cG_{\alpha_k}}\sum_{i\in A} |a_i|
\le\frac{16}{c^3} \max_{k\in\N,A\in\cG_{\alpha_k}} c^{3k}\sum_{i\in A} |a_i|
\end{equation}
(where the second inequality is an easy application of the inequality of Holder).

We now choose  (infinitely or finitely many) sequences $\xb^{(1)}$,   $\xb^{(2)}$ in $X$
 as follows.

We first choose $\xb^{(1)}\subset B_X$ weakly null, seminormalized,
 so that 
 $a(\xb^{(1)},\alpha_1)>c$. 
By Lemma \ref{L:6.3} (using also the
 part which starts with ``moreover'' for $\alpha_1$)  we can assume that the conclusions of
 Lemma \ref{L:6.3} are satisfied.

Now we consider the following two cases:

\noindent Case 1: For all $k\in\N$ it follows that $a(\xb^{(1)},\alpha_k)\ge c^{3k}$.
 In that case we are done for the moment and do not continue to choose other 
 sequences in $X$.

\noindent Case 2: There is a smallest $k_1$ (necessarily bigger than 1) for which
  $a(\xb^{(1)},\alpha_{k_1})< c^{3k_1}$.

In this case we choose a seminormalized and weakly null sequence (recall (\ref{E:6.4.1}))
\begin{equation}\label{E:6.4.4}
\yb^{(1)}\subset \frac{5}c c^{3k_1} B_X \text{ with } a(\yb^{(1)},\alpha_{k_1})\ge 5c^{3k_1}
\end{equation} 
and define $(x^{(2)}(n))=\xb^{(2)}=\xb^{(1)}+\yb^{(1)}$. 
Using now Lemma \ref{L:6.8} we find a subsequence $(x(2,n_i^{(2)}))_{i\in\N}$  of 
  $\xb^{(2)}=(x(2,n))_{n\in\N}$ so that
\begin{equation}\label{E:6.4.5}
a((x (2,n_i^{(2)}))_{i\in\N}, \alpha_{k_1})\ge c^{3k_1}.
\end{equation}
Once again we can consider two cases, namely

\noindent Case 1: For all $k\in\N$, $k\ge k_1$, it follows that $a((x(2,n_i^{(2)})),\alpha_k)\ge c^{3k}$.

\noindent Case 2: There is a smallest $k_2$ (must be  bigger than $k_1$)  for which
  $a(\xb(2,n^{(2)}_i),\alpha_{k_2})|\!<\! c^{3k_2}$.

As before we stop in Case 1 and define in Case 2 sequences $\yb^{(2)}$, $\xb^{(3)}$ and 
 $n^{(3)}_i$ as before.

We proceed in that way and eventually find an $\ell\in\N\cup\{\infty\}$,
   a strictly 
increasing sequence $(k_{j})_{1\le j<\ell+1}$   and 
 sequences $\yb^{(j)}\subset \frac{5}c c^{3k_j}B_X$, and 
$\xb^{(j)}=\xb{(1)}+\yb^{(1)}+\ldots \yb^{(j-1)}\subset X$, 
   and  subsequences of $N$
 $=\N=(n^{(1)}_i)\supset (n^{(2)}_i)\ldots$, for $1\le j<\ell+1$ so that for  any 
 $j\in\N$,  $1\le j<\ell+1$, it follows that
\begin{equation}\label{E:6.4.6a}
(\yb^{(j)})\subset\frac{5}c c^{3k_j} B_X
\end{equation}
satisfying the conclusions of Lemma \ref{L:6.3} and 
\begin{equation}\label{E:6.4.6}
 a((x(j,n^{(j)}_i)),\alpha_k)\ge c^{3k},
\text{ whenever } k_{j-1}\le k<k_j
\end{equation}
(where $k_0=1$ and $k_\ell=\infty$ if $\ell<\infty$). Moreover 
 (\ref{E:6.4.6}) holds for all subsequences of  $(x^{(j)},n^{(j)}_i))_{i\in\N}$.

Now we choose $\xb$ to be the diagonal sequence of $x(j,n^{(j)}_i)$, if $\ell=\infty$
(i.e. $x_i=(x(i,n^{(i)}_i))_{i\in\N}$)  or we choose
 $\xb= x(\ell,n^{(\ell)_i})$ if $\ell<\infty$.

For any $k\in\N$, there is a $j\in\N$, $1\le j<\ell+1$, so that
 $k_{j-1}\le k< k_j$ and we notice from (\ref{E:6.4.6a}) that
 $\xb$ is up to some finitely many elements  a subsequence 
 of $x^{(j-1)}$ and a sequence whose elements are of norm  not 
 exceeding $\frac{5}c\sum_{j\ge k_j}  c^{3j}$. 
 Thus, $\xb$ is  seminormalized and we 
  conclude from  (\ref{E:6.4.6}) that
 \begin{equation}\label{E:6.4.7}
a(\xb,\alpha_k)\ge c^{3k}-\frac{5}c\sum_{i\ge k_j}  c^{3i}\ge \frac12c^{3k},
\text{ whenever }j\in\N\text{ and }k_{j-1}\le k<k_j.
 \end{equation}

Using Lemma \ref{L:6.7} again we can assume  (by passing to subsequence) that
 for any $(a_i)\in\coo$  
 \begin{equation}\label{E:6.4.8}
 \Bigl\|\sum a_i x_i\Bigr\| \ge 
 \frac1{16}\max_{k\in\N,A\in\cG_{\alpha_k}}c^{3k}\sum_{i\in A} |a_i|. 
 \end{equation}
Combining (\ref{E:6.4.2}) and  (\ref{E:6.4.8})
  we showed that for any $(a_i)\in\coo$ it follows that 
\begin{equation*}
\frac{c^3}{16}  \Bigl\|\sum a_i y_i\Bigr\|
\le\max_{k\in\N,A\in\cG_{\alpha_k}} c^{3k}\sum_{i\in A} |a_i|
\le 16\Bigl\|\sum a_i x_i\Bigr\|.
\end{equation*} 
Thus the mapping $x_i\mapsto y_i$ extends
to a linear bounded operator $T$ which factors through a space which
 is hereditarily $c_0$  (by Lemma \ref{L:6.9}).  On the other hand we assumed that
$X$ does not contain $c_0$, therefore $T$ must be strictly singular. \end{proof}

\noindent Thomas Schlumprecht, Texas A\&M University, College Station, Tx 77843.
       
\noindent email: {\tt schlump@math.tamu.edu}
\end{document}